\documentclass{amsart}
\usepackage{amsmath,amsthm,amssymb,amscd, amstext, amsopn, amsxtra}

\renewcommand{\binom}[2]{\genfrac{(}{)}{0pt}{}{#1}{#2}}
\newtheorem{thm}{Theorem}[section]
\newtheorem{lemma}[thm]{Lemma}
\newtheorem{cor}[thm]{Corollary}
\newcommand{\C}{{\mathbb C}}
\newcommand{\Z}{{\mathbb Z}}
\newcommand{\N}{{\mathbb N}}
\newcommand{\F}{{\mathbb F}}

\newcommand{\Q}{{\mathbb Q}}
\newcommand{\D}{{\mathbb D}}
\newcommand{\1}{{\bf 1}_\lambda}

\newcommand{\U}{{\it U}}
\newcommand{\n}{{\eta}}
\newcommand{\B}{B_\lambda}
\newcommand{\Ca}{{\mathfrak C}}

\newcommand{\A}{{\mathfrak A}}

\newtheorem{prop}[thm]{Proposition}
\newtheorem{defn}{Definition}
\theoremstyle{remark}
\newtheorem{rem}{Remark}
\newtheorem{warning}{Warning}

\DeclareMathOperator{\Rep}{Rep}
\DeclareMathOperator{\Proj}{Proj}
\DeclareMathOperator{\Res}{Res}
\newcommand{\Hom}{{\mathop{\mathrm{Hom}}\nolimits}}

\newcommand{\Ext}{{\mathop{\mathrm{Ext}}\nolimits}}
\newcommand{\Ind}{{\mathop{\mathrm{Ind}}\nolimits}}
\newcommand{\aff}{{\mathop{\mathrm{aff}}\nolimits}}
\newcommand{\Indhat}{{\mathop{\widehat{\mathrm{Ind}}}\nolimits}}

\newcommand{\fin}{{\mathop{\mathrm{fin}}\nolimits}}

\newcommand{\limitm}{{\mathop{\varinjlim}\limits_m }\:}
\newcommand{\limitlm}{\mathop{\varprojlim}\limits_m \: }
\newcommand{\isom}{\simeq}
\newcommand{\surj}{\twoheadrightarrow}
\newcommand{\inj}{\hookrightarrow}

\renewcommand{\P}[1]{{\mathcal P}_{#1}}
\DeclareMathOperator{\ad}{ad}
\DeclareMathOperator{\tr}{tr}

\renewcommand{\Im}{\mathcal Im}
\newcommand{\cc}{{: }}
\newcommand{\ev}{{\mathop{\mathrm{ev}}\nolimits}}

\newcommand{\evst}{\ev^\ast}
\DeclareMathOperator{\pr}{pr}
\newcommand{\prlam}{\pr_\lambda}
\newcommand{\prhat}{\check{\pr}}
\newcommand{\e}{e}
\newcommand{\f}{f}
\newcommand{\h}{h}
\newcommand{\est}{\e^\ast}
\newcommand{\hst}{\h^*}
\newcommand{\hsti}{\h_i^*}
\newcommand{\esti}{{\e_i^\ast}}
\newcommand{\ed}{\e_i^{(n)}} 
\newcommand{\edst}{\e_i^{(n) \ast}}
\newcommand{\fsti}{{\f_i^\ast}}
\newcommand{\fst}{\f^\ast}
\newcommand{\fhat}{{\check \f}_i^\ast}
\newcommand{\fd}{\f_i^{(n)}}
\newcommand{\fdst}{\f_i^{(n) \ast}}
\newcommand{\et}{{\widetilde{\e}_i}} 
\newcommand{\ft}{\widetilde{\f}_i}
\newcommand{\ftd}{\ft^{(n)}}
\newcommand{\epsi}{\varepsilon_i}
\newcommand{\epsit}{\widetilde{\epsi}}
\newcommand{\nui}{\varphi_i}
\renewcommand{\nu}{\varphi}
\newcommand{\epsihat}{{\varepsilon_i^\wedge}}
\newcommand{\ethat}{{\et^\wedge}}
\newcommand{\fthat}{{\ft^\wedge}}
\newcommand{\esthat}{{{\esti}^\wedge}}

\DeclareMathOperator{\soc}{soc}
\DeclareMathOperator{\cosoc}{cosoc}
\DeclareMathOperator{\ch}{ch}

\newcommand{\RHaff}{\Rep H_n^{\aff}}
\newcommand{\Raff}{\Rep_q^{\aff}}
\newcommand{\RqHaff}{\Rep_q H_n^{\aff}}

\newcommand{\Rlam}{\Rep_q^{\lambda}}

\newcommand{\RHlam}{\Rep H_n^{\lambda}}

\newcommand{\Rfin}{\Rep_q^{\fin}}
\newcommand{\Haffn}{H_{n}^{\aff}}
\newcommand{\Haff}[1]{H_{#1}^{\aff}}
\newcommand{\Hlam}[1]{H_{#1}^{\lambda}}

\renewcommand{\sl}{\widehat{\mathfrak sl}_\ell}
\newcommand{\Rx}{R[X_1^{\pm 1}, \ldots, X_n^{\pm 1}] }

\newcommand{\triv}{{\bf 1}_{(1 q)}}
\newcommand{\stt}{{\mathop{{ \mathrm{\bf St}}}\nolimits}}
\newcommand{\st}[2]{\stt_{({#1} {#2})}}
\newcommand{\St}{\st{q }{ 1}}
\newcommand{\fq}[2]{\widetilde{\f}_{{#1}\: {#2}}}
\newcommand{\ech}[1]{{\varepsilon}_{#1}^\wedge}
\newcommand{\echone}{\ech{1}}
\newcommand{\wt}{{\mathop{\mathrm{wt}}\nolimits}}
\newcommand{\wti}{\wt_i}
\newcommand{\gr}{{\mathop{\mathrm{gr}}\nolimits}}
\newcommand{\Hdeg}{\overline{H}_n^{\gr}}
\newcommand{\Hldeg}{\overline{H}_n^{\lambda}}

\def\shuffle{\,\raise 1pt\hbox{$\scriptscriptstyle\cup{\mskip
               -4mu}\cup$}\,}

\def\<{\langle}
\def\>{\rangle}

\begin{document}
\title{Affine ${\widehat{\mathfrak sl}_p}$ controls the modular representation theory of the symmetric group and related Hecke algebras}
\author{I. Grojnowski}
\address{DPMMS, 16 Mill Lane, Cambridge, CB2 1SB}
\email{groj@dpmms.cam.ac.uk}

\abstract{
In this paper we prove theorems that describe how the representation
theory of the affine Hecke algebra of type $A$ and of related algebras
such as the group algebra $\F_lS_n$ of the symmetric group are
controlled by integrable highest weight representations of the 
{\it characteristic zero} affine
Lie algebra $\sl$. In particular we parameterise the representations
of these algebras by the nodes of the crystal graph, and give various
Hecke theoretic descriptions of the edges.

As a consequence we find for each prime $p$ a basis of the integrable
representations of $\sl$ which shares many of the remarkable
properties, such as positivity, of the global crystal basis/canonical
basis of Lusztig and Kashiwara.  This {\it $p$-canonical basis} is the
usual one when $p = 0$, and the crystal of the $p$-canonical basis is
always the usual one.

The paper is self-contained, and our techniques are elementary (no perverse
sheaves or algebraic geometry is invoked).
}
\endabstract

\maketitle



\section{Introduction}
\label{sec-intro}

This paper is concerned with the representation theory of the affine
Hecke algebra of type $A$ and of the cyclotomic Hecke algebras.

We parameterise the finite dimensional representations of these
algebras over any field $R$ and prove various results
about the behaviour of irreducible modules under restriction
and induction. Each of these algebras contains a large commutative
subalgebra, and we also describe how the failure of this algebra to
act semisimply controls the combinatorics of the representation
theory.

In contrast to the existing literature on these algebras, we prove
all our results without the use of sophisticated machinery or explicit
combinatorics---perverse sheaves and the geometry of the graded
nilpotent cones are notably absent from this work, as is the 
study of partitions.

In their place we prove the following theorem, which is a remarkable
rigidity property of the representation theory: there is an action of
the {\it characteristic zero}
 affine Lie algebra $\sl$ on the Grothendieck group of Hecke
algebra representations. Furthermore, the irreducible Hecke representations
define the natural crystal structure on the $\sl$ representation.

As an immediate consequence we recover the explicit combinatorics of
the Hecke algebra representation theory. For example, the simplest
case of these theorems identifies the Grothendieck group of symmetric
group representations in characteristic $p$ with an integral form of
the basic representation of $\widehat{\mathfrak sl}_p$. This
representation has a construction as a Fock space (the ``principal
realisation'').  The well known natural parametrisation by $p$-regular
partitions of the irreducible characteristic $p$ representations of
the symmetric group follows immediately.
In particular, this explains why the generating function for the number
of irreducible mod-$p$ representations of $S_n$ is just the character of 
the basic representation of $\widehat{\mathfrak sl}_p$.

To describe the results in more detail, consider the direct sum over
all $n$ of the Grothendieck group of representations of the affine
Hecke algebra of type $A_n$. (This, and all other terms, are carefully
defined in the body of the paper).

Similarly, consider the direct sum over all $n$ of the Grothendieck
group of representations of the symmetric group $S_n$.  Then it is a
classical observation that these are cocommutative Hopf algebras;
for the symmetric group this has been rediscovered many times (see \cite{Mac}),
but for the affine Hecke algebra this is due to \cite{BZ}.
In theorem \ref{uz-bi} we identify this algebra---it is just the dual
to the enveloping algebra $\U\n_\ell$ of the upper triangular part of 
the affine algebra $\sl$; here $l$ is the order of the parameter $q \in R^\times$
which enters in the definition of the affine Hecke algebra.

More generally, consider the cyclotomic Hecke algebra defined by Ariki
and Koike \cite{AK}.  This is a deformation of the group algebra of
the wreath product of the symmetric group $S_n$ with a cyclic group of
order $r$. The deformation depends on an $r$-tuple of elements of
$R^\times$, $\lambda = (q_1,\dots, q_r)$, and we denote the corresponding
algebra $\Hlam{n}$. When $r = 1$ and $\lambda = (1)$ this is
just the finite Hecke algebra.

If we now sum the Grothendieck groups of representations of $\Hlam{n}$
for fixed $\lambda$ it is no longer true that this is a Hopf algebra.
However it is obviously a comodule for the Hopf algebra (dual to
$\U\n_\ell$) built out of the affine Hecke algebra. The first of our
main theorems, theorem \ref{uz-hw}, says that it has many more
symmetries---that it is in fact dual to a module for the entire affine
algebra $\sl$, and moreover that this module is an irreducible
integrable highest weight module with highest weight determined by
$\lambda$.

Even in the classical case of the symmetric group and its deformations
($r =1$) this is new information: it identifies the Hopf algebra built
out of $S_n$ with the principal realisation of the basic
representation of $\sl$, and it identifies the action of $\sl$ by
vertex operators on this representation. For other cyclotomic Hecke algebras it 
extends the results of \cite{A}.

To prove this theorem we must introduce several new ingredients. The
first is the action of the Chevalley generators $f_i$ of the lower
triangular part $\n^-_\ell$ of $\sl$.  Unlike the operators $e_i$ of
$\n_\ell$ which have been known since the 1950s, and which arise in an
obvious way from the affine Hecke algebra, the definition of the $f_i$
is new to this paper. After a variety of preliminary results on the
affine Hecke algebra, we begin the study of these operators in section
\ref{sec-ef}.

The operators $e_i$ and $f_i$ are defined directly on the module
category, but will not satisfy the defining relations of $\sl$ before
we pass to the Grothendieck group. However, by considering the cosocle
filtration of these operators (i.e.~ before passing to the
Grothendieck group), we can define the {\it leading term} of the
operators $e_i$ and $f_i$. These leading terms have a beautiful
interpretation as ``crystal operators'' \cite{Ks}, and allow us to
define the {\it crystal graph} structure of the representations. This
structure generalises the classical ``branching laws'' for
representations of the symmetric group.

The crystal graph is a graph with nodes given by irreducible
representations, and with an edge between irreducible representation
$M$ of $\Hlam{n}$ and $N$ of $\Hlam{n-1}$ if $N$ occurs in the cosocle
of the restriction of $M$.  The operators $e_i^*$ which are dual to
$e_i$ refine restriction, and we can label this edge with an $i$ if
$N$ occurs in $e_i^* M$.

Our second main theorem (theorem \ref{thm-crystal}) is that this graph is
{\it precisely} the usual crystal graph of the representation determined by
$\lambda$. This gives a combinatorial parametrisation of
irreducible representations, and shows this parametrisation depends
only on $l$---the order of $q$ in $R^\times$.  (When $R =\C$ the
modules for $H_n^{\aff}$ were first parameterised in \cite{BZ} when $l
= \infty$, and in \cite{G} for arbitrary $l$. Subsequently Vigneras
conjectured \cite{V} that the parametrisation depends only on $l$, and not 
on $R$.)

\pagestyle{myheadings}
\markboth{I. Grojnowski}{Representations of Hecke algebras}

To prove this we must engage in a detailed study of the modules for the
affine Hecke algebra.  We begin by showing that $e_i^*M$ has
a simple cosocle.  This result, which is the main result of \cite{GV},
generalises the classical ``multiplicity one'' properties of
restriction for complex representations of the symmetric group, and
the corresponding property for characteristic $p$ representations
\cite{Kv}.


We then study the relationship between the crystal operators and the failure
of semisimplicity of both the Hecke algebra and a large commutative subalgebra of
it (the functions $\Rx$
on the maximal torus of $GL_n$). We find that though $e_i^*M$ is
a complicated module which is 
far from being semisimple, it has a uniserial
part to its composition series which admits a clean description in various ways.

In particular, the cosocle of $e_i^*M$ occurs in a uniserial chain
inside $e_i^*M$, and the length of the chain in which it occurs is
precisely the maximal size of a Jordan block for $X_n$ on $M$. 
We prove that this
length can also be read off the image of $M$ inside the Grothendieck
group of $\Rx$-modules.

These results, and the analogous (but harder) results for induction, form
the bulk of section \ref{sec-crystal}.

Finally, we can reinterpret the theorems \ref{thm-crystal}, \ref{uz-hw}
as defining a family of
interesting bases in the representations of $\sl$. Each of these bases
has the {\it same} crystal graph, but these {\it $p$-canonical bases}
are all different. In the case $R = \C$, the $0$-canonical basis
coincides with the basis defined by Lusztig and Kashiwara.  Each
$p$-canonical basis shares all of the remarkable properties of the
$0$-canonical basis---for example, the structure constants of $e_i$
and $f_i$ are non-negative integers, and the basis of the Verma
descends to a basis of the integrable representations.  These bases
are just the dual to the irreducible representations, and $p$ is the
characteristic of the field $R$. (The dual to the $p$-canonical basis
also has a Hecke theoretic interpretation---they are dual to the
projective representations).

To summarise, this paper ``explains'' all of the combinatorics of the
representation theory of the symmetric group and Hecke algebras---it is
just the combinatorics of the crystal graph of $\sl$. (Our proofs are
free of such combinatorics.)

\subsection{Acknowledgement} I would like to thank Monica Vazirani for 
her aid in preparing this paper. Her many suggestions have improved
both the mathematics and exposition, and our collaboration on
\cite{GV} has greatly increased my understanding of Hecke algebras and
their representation theory.

\section{Notation}
\label{sec-notate}

Throughout the paper, we fix a field $R$, and an invertible
element $q \in R$, i.e.,~a homomorphism $\Z[q, q^{-1}] \to R$.
We write $\mu_q = \{ q^i \mid i \in \Z \}$
for the powers of $q$ in $R$, and $\ell = | \mu_q |$,
so $\ell \in \N \cup \{ \infty \}$, and $\ell$ is the
order of $q \in R^\times$.

We further {\it assume that $q \neq 1$\/}. This is for clarity of
exposition; the changes in the statements of definitions and theorems
that must be made when $q = 1$ are sketched in section \ref{sec-q=1}.
The most interesting case when $q=1$ is the modular representation
theory of the symmetric group. Then $R= \F_p$ and $q=1$, but $\ell =
p$; see  section \ref{sec-q=1}.

Define the Dynkin diagram of $\mu_q$ to be the directed graph with
vertices the elements of $\mu_q$, and an edge $q^i \to q^j$ if
$q^{i-j-1} = 1$.  (A slightly classier notation, which we often use,
 is to write $i, j \in
\mu_q$ instead of $q^i, q^j \in \mu_q$, and then $i \to j$ if $i
j^{-1} = q$.)  The isomorphism type of this
graph depends only on $\ell$, and not on more general properties of
$q$ or $R$.  This feature will be mirrored by the properties of the
representation theory we study. (Conjecturally \cite{DJ} the
representation theory depends only on the characteristic of $R$ and the Dynkin
diagram $\mu_q$.)

The Dynkin diagram of $\mu_q$ defines an affine Lie algebra $\sl$.
The theorems of this paper are a description of how this
Lie algebra controls the representation theory of the affine
and cyclotomic Hecke algebras.  
We recall the definition and basic properties of 
$\sl$ and its representation theory (as found in
\cite{Kac}) in section \ref{sec-lie}. 

Most of the results about Hecke algebras hold for arbitrary
rings $R$ when appropriately formulated.  This is an 
easy exercise, but for clarity we have phrased results only for fields.

\subsection{Some common notation}
\label{sec-notation}

If $A$ is an $R$-algebra, we write 
$A$-mod for the category of {\it all\/} left $A$-modules, and 
$\Rep A$ for the category of left $A$-modules which are
finite dimensional as $R$-modules. Also write $\Proj A$ for the subcategory
of $\Rep A$ consisting of finite dimensional projective $A$-modules.
(If $R$ is an arbitrary ring, we would also need to define various
subcategories, such as the category of $A$-modules which are
projective as $R$-modules, and so on.)

We recall that the socle of a module $M$, denoted $\soc(M)$,
is the largest semisimple submodule of $M$, and that the cosocle of $M$,
denoted $\cosoc(M)$, is its largest semisimple quotient.

We write $S_n$ for the symmetric group on $n$ letters, 
$s_i = (i \hspace{2ex} i+1)$ for the simple transpositions, $\ell \cc  S_n \to \N$
for the length function.
\subsection{Grothendieck Group}
\label{sec-grothgroup}
If $\Ca$ is an abelian category, we write
$K(\Ca)$ for the Grothendieck group of $\Ca$.
This is the quotient of the free abelian group with
generators the objects $M \in \Ca$ by the ideal generated
by the elements $M_1 - M_2 +M_3$ for every short exact sequence
$$
0 \to M_1 \to M_2 \to M_3 \to 0
$$
in $\Ca$.  
If objects in $\Ca$ have finite length and unique 
composition factors, we write $K(\Ca)^{\ast}$ for
the topological dual of $K(\Ca)$; i.e.~for the linear
functions $f \cc  K(\Ca) \to \Z$ such that
$f(M) = 0$ for all but finitely many isomorphism classes
of irreducible objects $M \in \Ca$.

If $M \in \Ca$, let's write $[M]$ for its image in 
$K(\Ca)$.  Then as $M$ runs through the irreducible
objects in $\Ca$, the elements $[M]$ form a basis of
$K(\Ca)$, and the functions
$$
\delta_M\cc  K(\Ca) \to \Z \qquad \delta_M (N) =
\begin{cases} 0 & \text{if $M \not\isom N$}, \; \text{$N$ irreducible} \\
		1 & \text{if $M \isom N$} \end{cases}
$$
form a basis of $K(\Ca)^{\ast}$.
More generally, if $N \in \Ca$ and $M$ is an irreducible
object in $\Ca$, write $[M:N]$ for the multiplicity
of $M$ in a Jordan-Holder series of $N$, and extend this to
$[ \quad :\quad ]\cc  K(\Ca) \times K(\Ca) \to \Z$ by bilinearity.
Then write, for any $M \in  K(\Ca)$, $\delta_M \cc  K(\Ca) \to \Z $
for the function $N \mapsto [M:N]$.

Now, if $F\cc  \Ca \to \Ca^\prime$ is an exact functor of abelian
categories, we get an induced $\Z$-linear map
$F\cc  K(\Ca) \to  K(\Ca^\prime)$, and we can define its transpose
$F^{\ast}\cc  K(\Ca^\prime)^{\ast} \to  K(\Ca)^\ast$
by $F^{\ast} f = f F$.

We will apply all this to the category $\Rep A$ of finite
dimensional representations of an algebra $A$ (over $R$).
Suppose that $A$ and $A^\prime$ are two such algebras,
and that the cosocle of $A$ and $A'$ are direct sums of matrix algebras
over $R$; i.e.{} that they are separable algebras. Then the
the  irreducible representations of 
$A \otimes A^\prime$ are of the form 
$M \boxtimes M^\prime$, where $M$ is an irreducible $A$-module,
$M^\prime$ is an irreducible $A^\prime$-module.
More generally, we recall that under such assumptions
\begin{lemma}
$K(A \otimes A^\prime\text{-mod}) = K(A\text{-mod}) \otimes K(A^\prime\text{-mod})$
\label{lemma-tensor}
\end{lemma}
\noindent
which is certainly {\it not\/} true before passing to the
Grothendieck group.
(Here, the tensor product is algebraic---elements consist
of finite linear combinations of the elements 
$[M] \boxtimes [N]$ 
as $[M], [N]$  run through a basis
of $K(A\text{-mod})$ and $K(A^\prime\text{-mod})$ respectively.)

Write $K(\Ca)_\Q = K(\Ca) \otimes_\Z \Q$.
As $K(\Ca) $ is a torsion free $\Z$-module, 
$K(\Ca)_\Q$ is a $\Q$-vector space with 
distinguished sublattice $K(\Ca) \subset K(\Ca)_\Q$.

\section{Summary of properties of the affine Lie algebra $\sl$}
\label{sec-lie} 

First suppose $\ell \in \N$.
If $A$ is a ring, we write ${\mathfrak sl}_\ell (A)$ for
the Lie algebra of trace zero $\ell \times \ell$ matrices
over $A$ and $\sl(A)$ for the central extension of 
${\mathfrak sl}_\ell(A[t,t^{-1}])$ by $A$
$$ 0 \to A \cdot c \to \sl(A) \to {\mathfrak sl}_\ell(A[t,t^{-1}]) \to 0.
$$
This has Lie bracket $[f,g] = (fg - gf) + \tr \Res(
\frac{df}{dt} \cdot g) c$,
where $f, g \in {\mathfrak sl}_\ell(A[t,t^{-1}])$ and $\Res$
denotes the coefficient of $t^{-1}$.

If $\ell = \infty$, lets abuse notation and write 
$\sl(A)$ for the Lie algebra 
$${\mathfrak gl}_\infty(A)$$
of infinite matrices in which only finitely many entries
in any row or column are non-zero.
With this convention, the rest of this section is valid
for $\ell \in \N \cup \{ \infty \}$.

Write $\U_\Q \sl$ for the enveloping algebra of $\sl(\Q)$.
This has generators
$$ \e_i, \f_i, \h_i \qquad i \in \mu_q
$$
and relations
\begin{gather*}
[\e_i, \f_j] = \delta_{ij} \h_i \qquad 
[\h_i, \e_j] = c_{ij} \e_j        \qquad
[\h_i, \f_j] = -c_{ij} \f_j \\
(\ad \e_i)^{1 - c_{ij}} \e_j = 0 \qquad
(\ad \f_i)^{1 - c_{ij}} \f_j = 0   
     \end{gather*}
where $c_{ij} = 2 \delta_{ij} - (\delta_{ij^{-1}, q} + \delta_{ji^{-1}, q})$
is the Cartan matrix of $\sl$.

Write $\U_\Z \sl$ for Kostant's integral form of $\U_\Q \sl$.
This is a Hopf algebra over $\Z$, contained in $\U_\Q \sl$
as a lattice.
It is generated as an algebra by the elements
$$\ed = \frac{\e_i^n}{n!}, \quad \fd = \frac{\f_i^n}{n!},
\binom{\h_i}{n} = \frac{\h_i (\h_i -1) \cdots (\h_i -n +1)}{n!},
\qquad i \in \mu_q, n \in N
$$
with relations induced from the relations for $\e_i, \f_i$, and
$\h_i$.
(It is possible to write the relations explicitly; 
but we will not need this.)

We write $\U_\Z \n_\ell$ for the Hopf subalgebra of $\U_\Z \sl$
generated by $\ed$, $i \in \mu_q, n \in N$,
and $\U_\Z \n_\ell^-$  for the Hopf subalgebra of $\U_\Z \sl$
generated by $\fd$, $i \in \mu_q, n \in N$.

Recall that a representation $V$ of $\U_\Q \sl$
is called {\it integrable of lowest weight\/} if
\begin{enumerate}
\item[(i)] Each $\h_i$ acts semisimply
\item[(ii)] Each $\e_i$ and $\f_i$ acts locally nilpotently
\item[(iii)] There exists a finite dimensional subspace
$W \subseteq V$ such that $\U_\Q \n_\ell \cdot W = V$
\end{enumerate}
and that the category of such representations is semisimple.

Given an integrable lowest weight module $V$, the space
of invariants for $\U_\Z \n_\ell^-$ ({\it lowest weight
vectors\/}) is preserved by each $\h_i$.
If this space is one dimensional $V$ is irreducible, 
and if $\1$ is a lowest weight vector with 
$\h_i \1 = \lambda_i \cdot \1 \quad \forall i$,
then $\lambda_i \in \N$ and $\sum \lambda_i < \infty$.

This sets up a correspondence between 
$$(\text{functions } \lambda\cc  \mu_q \to \N, \sum \lambda_i < \infty)
\qquad \text{and} \qquad
\left( \begin{array}{l} \text{irreducible integrable} \\ \text{lowest weight modules.}
\end{array} \right)
$$
Given $\lambda$, write $L_\lambda$ for the corresponding
irreducible module.

As a $\U_\Q {\n}_\ell$-module, $L_\lambda$ is generated
by a lowest weight vector $\1$, and
$$L_\lambda = \left. \U_\Q {\n}_\ell \right/ {\U_\Q {\n}_\ell
\langle \e_i^{\lambda_i +1} \mid i \in \mu_q \rangle}.
$$
In particular, as the intersection of the left ideals
$\bigcap_{\lambda: \mu_q \to \N} \U_\Q {\n}_\ell 
\langle \e_i^{\lambda_i +1} \mid i \in \mu_q \rangle$
is trivial, it follows that
if $x \in \U_\Q {\n}_\ell$ acts as zero on every
integrable lowest weight representation,
then $x=0$.

We write $\Lambda_0\cc  \mu_q \to \N$ for the function
$\Lambda_0 (i) = \delta_{i,1}$.
The corresponding representation is called the 
{\it basic representation\/} of $\sl$.
More generally, define the fundamental weights $\Lambda_j : \mu_q \to \N$
by $\Lambda_j(i) = \delta_{i,j}$ and define the roots $\alpha_i : \mu_q \to \Z$
by $\alpha_i = 2\Lambda_i - \Lambda_{qi} - \Lambda_{q^{-1}i}$.

Each $L_\lambda$ carries a non-degenerate symmetric bilinear
form, the Shapovalov form $( \quad, \quad)$, which is determined by
requiring
$$(\1, \1) = 1$$
and 
$$(\e_i x, y) = (x, \f_i y)	\qquad \text{for all $i \in \mu_q$.}$$

\section{Definitions and first properties of Hecke algebras}
\label{sec-defns}

\subsection{}
\label{sec-hfin}
The {\it finite Hecke algebra\/}, $H_n$ is the $R$-algebra 
with generators
$$T_1, \ldots , T_{n-1}$$ 
and relations
\begin{eqnarray}
\label{braid}
 \text{{\it braid relations\/}} \hspace{1ex} &
 T_i T_{i+1} T_i = T_{i+1} T_i T_{i+1}, \quad
T_i T_j = T_j T_i,\; | i-j | > 1 &
\\
\label{quad}
\text{{\it  quadratic relations\/}  } &
(T_i + 1)(T_i - q) = 0. & 
\end{eqnarray}

The braid relations imply that if 
$w = s_{i_1} \cdots s_{i_r}$ and $\ell (w) = r$,
then $T_{i_1} \cdots T_{i_r}$ depends only on $w \in S_n$.
It is denoted $T_w$, and the $T_w, w \in S_n$ form a basis
of $H_n$ over $R$.

The  {\it affine Hecke algebra\/} $H_n^{\aff}$ is the $R$-algebra,
which as an $R$-module is isomorphic to
$$H_n \otimes_R \Rx.$$
The algebra structure is given by requiring that $H_n$ and 
$R [ X_i^{\pm 1} ]$ are subalgebras, and that
\begin{equation}
T_i X_i T_i = q X_{i+1} . 
\end{equation}
Equivalently, if $f \in \Rx$, then
\begin{equation}
\label{Tf-fT}
T_i f -  {}^{s_i}\!f T_i = (q-1)\frac{f - {}^{s_i}\!f}  {1- X_i X_{i+1}^{-1}}
\end{equation}
where $s_i \in S_n$ acts on $\Rx$
by permuting $X_i$ and $X_{i+1}$.

The affine Hecke algebra so defined is an extension by a Laurent
polynomial algebra of the Hecke algebra associated to the Coxeter
group with Dynkin diagram $\mu_q$.  This definition and the isomorphism
is due to Bernstein.
\begin{prop}{(Bernstein)} The center of $H_n^{\aff}$,
$Z(H_n^{\aff})$, is isomorphic to symmetric Laurent polynomials.
$$Z(H_n^{\aff}) \isom \Rx^{S_n}$$
\label{prop-center}
\end{prop}
\begin{proof}
The relation \eqref{Tf-fT} makes it clear that $\Rx^{S_n}$
is contained in the center.
Conversely, suppose $h = \sum T_w f_w \in Z(H_n^{\aff})$ where $f_w \in \Rx$.
Let $u \in S_n$ be maximal with respect to Bruhat order such that
$f_u \neq 0$.  If $u \neq 1$ then there exists some $i$ such that
$u(i) \neq i$.  As $X_i T_w = T_w X_{w^{-1}(i)} + \sum_{w' < w} T_{w'} g_{w'}$
for some $g_{w'} \in \Rx$,
the coefficient of $T_u$ is different in
 $X_i h$ and $h X_i$.
Hence $h = f_1 \in \Rx$.  But then \eqref{Tf-fT} implies $h$ is $S_n$-invariant.
\end{proof}

\begin{cor}
If $M \in H_n^{\aff}$-mod is absolutely irreducible, then $M$ is finite
dimensional, and in fact $\dim_R M \le n !$.
\label{cor-dim}
\end{cor}
\begin{proof}
As $\Rx$ is a free $\Rx^{S_n}$-module of rank $n !$,
$H_n^{\aff}$ is a free module over $Z(H_n^{\aff})$ of rank 
$(n !)^2$.  Dixmier's version of Schur's lemma implies that the center
of $H_n^{\aff}$ acts by scalars on absolutely irreducible modules, and hence
$M$ is an irreducible  module for a finite 
dimensional algebra of dimension $(n !)^2$.
\end{proof}
Suppose $R$ is algebraically closed. 
Then
the characters (i.e.~one dimensional representations) of the
center $Z(H_n^{\aff})$ are the orbits of $S_n$ on 
$(R^\times)^n =$  $ \text{Spec} R[X_1^{\pm 1}, \ldots,\linebreak[3]
 X_n^{\pm 1}]$.
Given any finite dimensional module $M \in \RHaff$, we can write  $M$
as a direct sum of {\it generalized\/} eigenspaces for the commutative
subalgebra $\Rx$, say
$$M = \bigoplus_{s \in (R^\times)^n} M [s]$$
where $s = (s_1, \ldots, s_n)$ and 
 $M[s] = \{m \in M \mid (X_i - s_i)^{\dim M} m = 0, 1 \le i \le n \}.$

If we also write
$$M_s = \{m \in M \mid X_i m = s_i m, 1 \le i \le n \}$$
for the  {\it actual\/}  simultaneous eigenspace of the elements
$X_i$, we have
$$M_s \neq 0  \iff M[s] \neq 0.$$
Further, if $M$ is irreducible, then 
$Z(H_n^{\aff})$ acts by the central character $S_n \cdot s
\in (R^\times)^n / S_n$ if and only if there exists some 
$w \in S_n$ such that $M_{ws} \neq 0$.

To generalize this to all of $\RHaff$, recall
the following general property of algebras.
\begin{lemma}
Let $A$ be an $R$-algebra, and $Z \subset A$ a central subalgebra.
If $M,  N$ are two $A$-modules on which $Z$ acts by  {\it different\/}
one dimensional characters, then 
$$\Ext_A^i(M,N) = 0,\qquad \text{for all i.}$$
\label{lemma-exts}
\end{lemma}
Now, let $M$ be an indecomposable finite dimensional $H_n^{\aff}$-module.
This has a finite filtration by irreducible $H_n^{\aff}$-modules,
and the lemma implies $Z(H_n^{\aff})$ acts by the same character
on each subquotient.  It follows that

\begin{prop}
\label{prop-blocks}
$\RHaff \isom \bigoplus_{s \in (R^\times)^n / S_n} \Rep_s H_n^{\aff}$
(direct sum of categories). A module $M$ is in $\Rep_s H_n^{\aff}$ if and
only if the support of $M$ as a $Z(H_n^{\aff})$-module
is $s \in (R^\times)^n/S_n$, i.e.~%
if and only if 
(i) there exists an $s^\prime \in (R^\times)^n$  with
$M_{s^\prime} \neq 0$, and $S_n \cdot s^\prime = s$,
and (ii)
if $M_{s^\prime} \neq 0$, then $S_n \cdot s^\prime = s$.
\end{prop}
The summands above are called  {\it blocks\/} of the category $\RHaff$.
If $s \in (R^\times)^n$, $\Rep_s H_n^{\aff} = 
\varinjlim_k \Rep (H_n^{\aff} / Z_s^k H_n^{\aff})$,
where $Z_s = \{ f \in \Rx^{S_n} \mid f(s) = 0\}.$ If $N \in \Rep_s H_n^{\aff}$
we say that $N$ has {\it central character} $s$.

\subsection{}
\label{sec-cyc}
Fix a function $\lambda \cc  \mu_q \to \Z_+$ such that $\sum_{i \ge 0} 
\lambda_i = r < \infty$.

The  {\it Ariki-Koike algebra\/}, or  {\it cyclotomic Hecke algebra\/},
is the $R$-algebra with generators
$$T_1, \ldots, T_{n-1} \qquad  \text{ and } \qquad T_0$$
and relations
$$\prod_{q^i \in \mu_q} (T_0 - q^i)^{\lambda_i} = 0$$
$$T_0 T_1 T_0 T_1 = T_1 T_0 T_1 T_0$$
as well as the braid relations \eqref{braid} and quadratic
relations \eqref{quad} in the definition of the 
finite Hecke algebra above for $i \ge 1$.

In particular, if $r =1$ then this is just the finite
Hecke algebra, and if $r =2$ this is the Hecke algebra
of type $B_n$ or $C_n$ with possibly unequal parameters.
The finite Hecke algebra is always a subalgebra of
the cyclotomic Hecke algebra.

There is a surjective algebra homomorphism,
first defined by Cherednik
$$\ev = \ev_{\lambda} \cc H_n^{\aff} \to H_n^\lambda$$
defined on the generators by $T_i \mapsto T_i, \qquad 1\le i\le n-1$
$$X_1 \mapsto T_0, \qquad  X_i \mapsto
q^{1-i} T_{i-1} \cdots T_1 T_0 T_1 \cdots T_{i-1} $$
Write $\evst \cc  \RHlam \to \RHaff$ for the induced map of
modules.

The image $\ev(\Rx)$ form a commutative subalgebra of 
$H_n^\lambda$, (``Murphy operators'') and
as $\ev$ is surjective, the image $\ev(Z(H_n^{\aff}) )$ is 
contained in the center of $H_n^\lambda$.
This implies the category of $H_n^\lambda$-modules splits
up into a direct sum, indexed by characters of $Z(H_n^{\aff})$,
$$\RHlam = \bigoplus \Rep_s H_n^{\lambda}.$$
The key result about these algebras is:
\begin{prop}[Ariki-Koike]
\label{prop-cycl}
The algebra $H_n^\lambda$ is finite dimensional, of dimension
$r^n n !$, where $r = \sum_{i \in \mu_q} \lambda_i$.
\end{prop}

The image under $\ev$ of the elements
$$ X_1^{a_1} \cdots X_n^{a_n} T_w, \qquad 0 \leq a_i <  r, w \in S_n$$
form a basis of $H_n^\lambda$.


Consider the modules $M$ for $H_n^{\aff}$ such that the only 
eigenvalues of $X_1$ on $M$ are powers of $q$.
Such modules form a full subcategory of $\RHaff$ which we denote
$\RqHaff$.

If $0 \to M_1 \to M_2 \to M_3 \to 0$ is exact in $\RHaff$,
and any two of the modules $M_i$ are in $\RqHaff$,
then so is the third, i.e.~$\RqHaff$ is closed under subquotients
and extensions.

Further, $\RqHaff$ is a direct sum of blocks of the category
$\RHaff$; i.e.

\begin{lemma}
\label{lemma-repqblocks}
If $M \in \RqHaff$, and $N \in \RHaff$ with $\Ext^\bullet (M,N) \neq 0$
or $\Ext^\bullet (N,M) \neq 0$, then $N \in \RqHaff$ also.
More precisely,
$$\RqHaff = \bigoplus_{s \in {\mu_q}^n /S_n} \Rep_s H_n^{\aff}.$$
\end{lemma}
This is immediate from the following lemma, and the above
description of the center of $H_n^{\aff}$.
\begin{lemma}
\label{lemma-Xq} 
If $M$ is a module for $H_n^{\aff}$ and there is some $i$ such that the only
eigenvalues for $X_i$ on $M$ are powers of $q$, then for all $j$,
the only eigenvalues of $X_j$ on $M$ are powers of $q$.
\end{lemma}
\begin{proof}
Write $X=X_{i-1}, T=T_{i-1}$.
Let $v$ be an eigenvector for $X$ and  for
$TX T = qX_{i}$. It is enough to show that the eigenvalues of $X$ and $TXT$ on
$v$ differ by a power of $q$. Consider the space spanned by $v$ and $Tv$. This is $X$-stable. If it is two dimensional, then with respect to the basis $v,Tv$ we have $T$ has matrix
$ \left( \matrix 0 & q \\  1 & q-1 \endmatrix \right) $ and
$X$ has matrix
$ \left( \matrix \mu & a \\  0 & \mu' \endmatrix \right) $.
It follows that  $TX T$ has matrix
$ \left( \matrix q\mu' & q\mu'(q-1) \\  \delta & q\mu + (q-1)\delta \endmatrix \right) $
where $\delta = a + \mu'(q-1)$. 
By assumption, $\delta = 0$, so $TX T$ has eigenvalues $q\mu, q\mu'$
and we are done. If $Tv $ is a multiple of $v$, then either $Tv= -v$ or $Tv=qv$.
In the first case $TX T$ has eigenvalue $\mu$, in the second case $q^2\mu$, and we're done.
\end{proof}

Given $\lambda\cc \mu_q \to \Z_+$, $\sum \lambda_i < \infty$,
it is immediate from the definition of $\Hlam{n}$ that if 
$M \in \RHlam,$  $\; \evst M \in \RqHaff$.
Define $\Rep_q^\lambda H_n^{\aff}$ to be the full subcategory of
$\RHaff$ whose objects are the modules $M$ such
that the Jordan blocks of $X_1$ on $M$ with eigenvalue $q^i$
have size less  than or equal to 
$\lambda_i$ (and there
are no other eigenvalues), i.e.~$\Rep_q^\lambda H_n^{\aff}$
consists of modules annihilated by $\prod (X_1 - q^i)^{\lambda_i}$.

Then if $M \in \Rep_q^\lambda H_n^{\aff}$, and $N$ is a subquotient
of $M$, then $N \in \Rep_q^\lambda H_n^{\aff}$ also.
However, extensions of modules in $\Rep_q^\lambda$ need
not be in $\Rep_q^\lambda$.

We clearly have
$$\evst \cc  \RHlam \to \Rep_q^\lambda H_n^{\aff}$$
is an equivalence of categories.

Given any $M \in \RqHaff$, there are infinitely many $\lambda$
such that $M \in \Rep_q^\lambda H_n^{\aff}$. 
More precisely, define a partial order on $\lambda$ by 
$\lambda \ge \mu$ if $\lambda_i \ge \mu_i$ for all $i$. Then
$$\mu \le \lambda \implies \Rep_q^\mu H_n^{\aff}
\subset \Rep_q^\lambda H_n^{\aff}$$
and
$$\RqHaff = \varinjlim_\lambda 
\Rep_q^\lambda H_n^{\aff}.$$

The exact functor $\evst = \evst_\lambda \cc 
\RHlam \to \RHaff$ has a left adjoint $\pr_\lambda$, and a right
adjoint $\prhat_\lambda$.
To define them, write
$$I_\lambda = \ker (\ev \cc  H_n^{\aff} \to H_n^\lambda).$$
As $\ev$ is an algebra homomorphism, $I_\lambda$ is a
two sided ideal, and as $\ev$ is surjective
$$H_n^\lambda = \left. H_n^{\aff} \right/ I_\lambda.$$
Now if $N \in \RHaff$, we have $\left. N\right/ I_\lambda N$ is an
$H_n^{\aff}$-module on which $I_\lambda$ acts trivially
and so $\left. N \right/ I_\lambda N$ is an $H_n^\lambda$-module.
If $M \in \RHlam$ we have
$$\Hom_{H_n^\lambda}(\left.N\right/ I_\lambda N, M) =
\Hom_{H_n^{\aff}}(\left.N\right/ I_\lambda N,\evst M) =
\Hom_{H_n^{\aff}}(N, \evst M) $$
and so if we define 
$$ \pr_\lambda (N) =\left. N\right/ I_\lambda N$$
and
$$\prhat_\lambda (N) = N^{I_\lambda}$$
we have proved that $\pr_\lambda$ is the left adjoint to 
$\evst$.  Similarly, $\prhat_\lambda$ is the right adjoint,
and neither functor is exact.

Define
$$\Raff = \bigoplus_{n \ge 0} \RqHaff$$
$$\Rlam = \bigoplus_{n \ge 0} \RHlam  \xrightarrow{\substack{\evst \\ \isom}} 
\bigoplus \Rep_q^\lambda H_n^{\aff}$$
$$\Rfin = \bigoplus_{n \ge 0} \Rfin \xrightarrow{\substack{\evst \\ \isom} } 
\bigoplus \Rep_q^{\Lambda_0} H_n^{\aff}.$$
We will investigate some rigid structures on these categories.

\section{Generalities on Induction and Restriction}
\label{sec-indres}

The results in this section are easy and are mostly well known
(with the possible exception of \ref{thm-enhanced}). We omit proofs which
are trivial variants of what is already in the literature.

\subsection{}
Recall that if $A \subset B$ are $R$-algebras, the exact functor
of restriction
$$\Res_A^B \cc  {B}\text{-mod} \to {A}\text{-mod} $$
has left and right adjoints, $\Ind$ and $\Indhat$ defined by
$$\Ind_A^B \cc  {A}\text{-mod} \to {B}\text{-mod}   \hspace{3ex} 
M \mapsto B \otimes_A M$$
$$\Indhat_A^B \cc  {A}\text{-mod} \to {B}\text{-mod}    \hspace{3ex}
 M \mapsto \Hom_A(B, M);$$
i.e. we have the Frobenius reciprocity
$$\Hom_B(\Ind_A^B M, N) = \Hom_A (M, \Res N) \qquad
\Hom_B(N, \Indhat M) = \Hom_A(\Res N, M).$$
If $B$ is a free $A$-module, then $\Ind$ and $\Indhat$
are exact functors also.
Further, if $A \subset B \subset C$ are inclusions of $R$-algebras,
we have transitivity of induction and restriction:
$$\Res_A^B \Res_B^C = \Res_A^C, \hspace{3ex} 
\Ind_B^C \Ind_A^B = \Ind_A^C,    \hspace{3ex}
\Indhat_B^C  \Indhat_A^B = \Indhat_A^C$$
Now apply these remarks to affine and cyclotomic Hecke algebras.
Given a sequence $P = (a_1, \ldots, a_k)$ of non-negative
integers, we have an obvious embedding 
$$H_{a_1}^{\aff} \otimes \cdots \otimes H_{a_k}^{\aff}
\hookrightarrow H_{a_1 + \cdots + a_k}^{\aff}$$
which makes
$H_{a_1 + \cdots + a_k}^{\aff}$ a free
 $H_{a_1}^{\aff} \otimes \cdots \otimes H_{a_k}^{\aff}$-module.
Applying the previous remarks we get exact functors 
$\Res, \Ind, \Indhat$.
These functors depend on the order $(a_1, \ldots, a_k)$ and not
just on the underlying set!

Write $\Rep_q(\Haff{a_1} \otimes \cdots \otimes \Haff{a_k})$
for the full subcategory of
 modules for $\Haff{a_1} \otimes \cdots \otimes \Haff{a_k}$
on which each $X_j$ acts with eigenvalues in $\mu_q$. 
Then the following is evident.
\begin{lemma}
\begin{enumerate}
\item[(i)]
 $\Res$ and $\Ind$ define functors $\Rep_q(\Haff{a_1} \otimes \cdots \otimes \Haff{a_k})
\leftrightarrows \Rep_q \Haff{a_1 + \cdots + a_k}$.
\item[(ii)]
 $K(\Rep_q(\Haff{a_1} \otimes \cdots \otimes \Haff{a_k}))
= K(\Rep_q \Haff{a_1}) \otimes \cdots \otimes K(\Rep_q \Haff{a_k})$.
\end{enumerate}
\end{lemma}
Similarly, we may define induction and restriction for finite
Hecke algebras, i.e.~functors
$$\Rep H_{P} 
\mathop{ \rightleftarrows}\limits_{\Res}^{\Ind}
 \Rep H_n.$$
It is clear that 
\begin{lemma}
$\evst \Res_{H_P}^{H_n} = \Res_{\Haff{P}}^{\Haffn} \evst$.
\end{lemma}
But note that though $\Ind_{\Haff{P}}^{\Haffn} \evst$ and
$\evst \Ind_{H_P}^{H_n}$ seem quite different, they are 
related. The results of section 
\ref{sec-detailedcrystal} describe the relation.

\subsection{Mackey formula}
\label{sec-mackey}
We now state the Mackey formula.  Unfortunately, it requires some
notation.

Given a sequence $P = (a_1, \ldots, a_k)$ of positive integers,
as above, with $\sum a_i = n$, write 
$\Haff{P} = \Haff{a_1} \otimes \cdots \otimes \Haff{a_k} \hookrightarrow
\Haffn$, $S_P = S_{a_1} \times \cdots \times S_{a_k} \hookrightarrow S_n$,
and $H_P = H_{a_1} \otimes \cdots \otimes H_{a_k}$.
Then if $P^{\prime}$ is another such sequence,
we define a  partial order on $S_P \left\backslash S_n \right/ S_{P^{\prime}}$ 
by setting $x \leq y$ if the maximal length coset representative of
$x$ is less than that of $y$ in Bruhat order.

Now if $x$ is a minimal length coset representative of an element 
$x \in S_P \left\backslash S_n \right/ S_{P^{\prime}}$, write $P \cap {}^x \!P^{\prime}$
for the new sequence of positive integers defined by listing
the parts of the  ordered partition of $1,\dots, n$ defined by 
$ P \cap {}^x \!P^{\prime}$.
Note $\Haff{P \cap {}^x \!P^{\prime}}$ is isomorphic to 
$\Haff{P^{\prime}\cap {}^{x^{-1}} \! P}$ by the isomorphism 
``conjugation by $x$'' which sends
$$T_w \mapsto T_{xwx^{-1}} \qquad	X_i \mapsto X_{x^{-1}(i)}.$$ 
\begin{lemma}
\begin{enumerate}
\item[(i)] $\Haffn$ is a free right $\Haff{P}$-module, and 
$$
\Haffn = \bigoplus_{w \in S_P \left\backslash S_n\right.} \Haff{P} \cdot T_w
$$
Similarly for $H_n$.
\item[(ii)]  If $x \in S_P \left\backslash S_n \right/ S_{P^{\prime}}$, then
$$
\Haff{P} \cdot T_x \cdot \Haff{P^{\prime}} =
\sum_{ \substack{y \le x \\ y \in S_P \left\backslash S_n \right/ S_{P^{\prime}}}}
\Haff{P} \cdot T_y \cdot \Haff{P^{\prime}}
$$
$$= \bigoplus_{\substack{a \le x \\ a \in S_P \left\backslash S_n\right.}}
\Haff{P} \cdot T_a
 = \bigoplus_{\substack{b \le x \\ a \in\left.  S_n \right/ S_{P^{\prime}} }}
T_b \cdot \Haff{P^{\prime}}.
$$
(The first isomorphism is as $\Haff{P} \text{-} \Haff{P^{\prime}}$-modules,
the second as $\Haff{P}$-modules, the last as
$(\Haff{P^{\prime}})^{\text{opp}}$-modules.)
\item[(iii)]  The above defines a filtration of $\Haffn$ by $(\Haff{P}, \Haff{P^{\prime}})$-bimodules.
The associated graded module
$$
\gr \Haffn = \bigoplus_{x \in S_P \left\backslash S_n \right/ S_{P^{\prime}}}
H_{P} S T_x H_{P^{\prime}}
$$
has summands isomorphic to  $\Haff{P} \otimes_{\Haff{P \cap {}^x \!P^{\prime}}} \Haff{P^{\prime}}$
as an $(\Haff{P}, \Haff{P^{\prime}})$-bimodule.
\end{enumerate} 
\end{lemma}
We omit the proof. It has the consequence
\begin{cor}
$\Res_{\Haff{P}}^{\Haffn} \Ind_{\Haff{P^{\prime}} }^{\Haffn} M$
admits a filtration with subquotients isomorphic to 
$$
\Ind_{\Haff{P \cap {}^x \!P^{\prime}}}^{\Haff{P}} w^{-1}
 \Res_{\Haff{P^{\prime}\cap {}^{x^{-1}} \! P}}^{\Haff{P^{\prime}}} M.
$$
\end{cor}
In the above filtration, $\Ind_{\Haff{P \cap P^{\prime}}}^{\Haff{P}} 
 \Res_{\Haff{P^{\prime}\cap  P}}^{\Haff{P^{\prime}}} M$
always sits as a subobject of $\Res \Ind M$.

If we apply the above corollary to $P= (1, \ldots, 1)$, we get a particularly
nice consequence.
Write 
$$\ch (M) = \sum_{s \in (R^\times)^n} \dim M[s] \cdot [s]
\in K(\Rep \Rx)$$
for the character of $M$ as a module for $\Rx$.
\begin{lemma}[shuffle lemma]
\label{shufflelemma}
If $M \in \RqHaff$, $N \in \Rep_q \Haff{m}$, then
$$
\ch \Ind(M \boxtimes N) = \sum_{s^{\prime}} (\dim M[s^{\prime}] \cdot
\dim N[s^{\prime \prime}])[s]
$$
where if $s = (s_1, \ldots , s_{n+m})$, $s^{\prime}$ is a 
subsequence $(s_{i_1}, \ldots, s_{i_n})$ where
$i_1 < \cdots < i_n$, and $s^{\prime \prime}$ is the
sequence obtained from $s$ by removing $s^{\prime}$.
In other words, the spectrum of $\Ind(M \boxtimes N)$
is obtained by {\it shuffling\/} the spectrum of $M$ and $N$.
\end{lemma}

\subsection{Boring central characters}
\label{sec-centralchars}

Let $P = (a_1, \ldots, a_k)$ be a sequence of positive integers,
$\sum a_i = n$, and write 
$\Haff{P} = \Haff{a_1} \otimes \cdots \otimes  \Haff{a_k} 
\hookrightarrow \Haff{n}$
as before. 
Write $S = \Rx$, so that $Z(\Haff{P}) = S^{S_P}$, and form the 
{\it enhanced\/} Hecke algebra
$$
\widetilde{\Haff{P}} = \Haff{P} \otimes_{Z(\Haff{P})} S.
$$
Clearly $Z(\widetilde{\Haff{P}}) = S$.
Define the $q$-discriminant,
$\Delta_q^P \cc  (R^\times)^n \to R$ by
$$
\Delta_q^P (s_1, \ldots, s_n) = \prod (s_i - q s_j)
$$
where the product runs over all pairs $(i,j)$ such that 
$1 \le i, j \le n$, and it is not the case that
both $i$ and $j$ lie in an interval 
$a_1 + \cdots + a_r +1, \ldots, a_1 + \cdots +a_r +a_{r+1}$.
Also write
$$
\Delta (s_1, \ldots, s_n) = \Delta_1^{(1,1,\ldots,1)} (s_1, \ldots, s_n)
= \prod_{i \neq j} (s_i - s_j)
$$
for {\it the\/} discriminant.
So if $\Delta (s) \neq 0$, there are $n !$ points in the 
$S_n$-orbit of $s$.
Now, if $s \in (R^\times)^n$, let's write 
$\Rep_s \widetilde{\Haff{P}}$ for the category of finite
dimensional $\widetilde{\Haff{P}}$-modules $M$ such that
the support of $M$ as an $Z(\widetilde{\Haff{P}}) = S$ module is
$s \in (R^\times)^n$.
\begin{thm}
\label{thm-enhanced}
\begin{enumerate}
\item[(i)] 
If $\Delta_q (s) \neq 0$, induction defines an equivalence
of categories
$$
\Rep_s \widetilde{\Haff{P}} \to \Rep_s \widetilde{\Haff{n}},
\qquad M \mapsto \Haff{n} \otimes_{\Haff{P}} M.
$$
\item[(ii)] 
If $\Delta (s) \neq 0$, there is an equivalence of categories
$$
\Rep_s \widetilde{\Haff{n}} \xrightarrow{\sim} \Rep_s \Haff{n}.
$$
\item[(iii)] 
Regardless, there is always an isomorphism
$$
K(\Rep_s \widetilde{\Haff{n}})  \xrightarrow{\sim} K(\Rep_s \Haff{n}).
$$
\end{enumerate}
\end{thm}
In particular, if $M$ is an irreducible $\Haff{n}$-module
with central character $S_n \cdot s$, then by defining
that $S$ acts via evaluation at $s$, we get an irreducible
$\widetilde{\Haff{n}}$-module, and conversely every irreducible
module in $\Rep_s \widetilde{\Haff{n}}$ is of such a form.

The following particular case will be of great importance to us.
\begin{cor}
\label{cor-kato}  \cite{Ka}
The $\Haff{n}$-module
$\Ind_{\Haff{1} \otimes \cdots \otimes \Haff{1}}^{\Haff{n}}
(q^i J_1 \boxtimes \cdots \boxtimes q^i J_1)$ is irreducible.
It is the unique irreducible $\Haff{n}$-module with central
character $(q^i \cdots q^i)$.
\end{cor}
Notice that the theorem tells us completely how to reconstruct
all of $\RHaff$ once we understand $\RqHaff$.
(This is one good reason for concentrating on $\RqHaff$
for the rest of the paper!)

\subsection{A useful criterion for irreducibility}
\label{sec-useful}
\begin{prop}
$\Ind_{\Haff{a} \otimes \Haff{b}}^{\Haff{a+b}}(M\boxtimes N)
\isom \Indhat_{\Haff{b} \otimes \Haff{a}}^{\Haff{a+b}}(N\boxtimes M)$.
\end{prop}
The preceding proposition, when combined with the following
useful observation, allows us to detect whether certain induced modules
are irreducible.
\begin{lemma}
\label{useful}
If $\Ind(M\boxtimes N) \isom \Indhat(M\boxtimes N)$, and
$M\boxtimes N$ is an irreducible module which occurs with 
multiplicity one in a composition series for 
$\Res \Ind(M\boxtimes N)$, then $\Ind(M\boxtimes N)$ is irreducible.
\end{lemma}
\begin{proof}
Let $0 \to K \to \Ind(M\boxtimes N) \to Q \to 0$ be an
exact sequence.
Then $M\boxtimes N$ is a submodule of $Q$, if $Q \neq 0$, by
Frobenius reciprocity, and as $\Ind(M\boxtimes N) \isom
\Indhat(M\boxtimes N)$, $M\boxtimes N$ is also a submodule of
$K$ if $K \neq 0$. As $M\boxtimes N$  occurs with multiplicity
one, either $K$ or $Q$ is zero, and so $\Ind(M\boxtimes N)$
is irreducible.
\end{proof}

\section{Examples}
\label{sec-examples}

The computations in this section will be used in sections \ref{sec-serre}
and \ref{sec-groth}.

\subsection{{ $H_1^{\aff}$}}
\label{sec-H1}
Recall $H_1^{\aff} = R [x, x^{-1}]$.

For $q^i \in \mu_q$, and $n \ge 1$, let's write
$q^i J_n$ for the rank $n$ Jordan block with 
eigenvalue $q^i$, i.e.
$$q^i J_n := R[x] / (x -q^i)^n.$$
This is an indecomposable $H_1^{\aff}$-module, and conversely
every indecomposable module in $\Rep_q H_1^{\aff}$ 
is of the form $q^i J_n$ for a unique $q^i \in \mu_q$ and $n \ge 1$.
Further 
$$\Hom(q^i J_n, q^j J_{n'}) = \Ext^1(q^i J_n, q^j J_{n^\prime}) = 
\begin{cases}
0 & \text{if $q^i \neq q^j$} \\
R^{min(n,n^\prime)} & \text{if $q^i = q^j$}
\end{cases}$$
and 
$$\Ext^k(q^i J_n, q^j J_{n^\prime}) = 0 \qquad \text{if $k \neq 0,1.$}$$
Hence $\Rep_q H_1^{\aff}$ is the direct sum of $\mu_q$-copies
of the same category, the category of 
 {\it finite Jordan blocks\/} (with a fixed eigenvalue).

\subsection{$\Haff{2}$-modules}
\label{sec-H2}
We list all the indecomposable $\Haff{2}$-modules.
First, an easy computation.
Write $S = R[X_1^{\pm1}, X_2^{\pm1}]$.
\begin{lemma}
\label{lem-h2}
Set $p = (1- X_1 X_2^{-1}) T_1 - (q-1)$. Then
\begin{enumerate}
\item[(i)] \label{lem-h2-i} $X_1 p = p X_2, \quad X_2 p = p X_1$
\item[(ii)] $p^2 = 
(q - X_1 X_2^{-1}) (q - X_2 X_1^{-1})$ is central.
\end{enumerate}
\end{lemma}
We describe $\Rep_{s_1, s_2} \Haff{2}$.
Choose $s = (s_1, s_2)$ in the orbit $\{ (s_1, s_2),
(s_2, s_1) \}$, and write $s' = (s_2, s_1)$.
Let $M$ be a module, and write $M = M_s + M_{s'}$ for
its decomposition into generalized eigenspaces for $S$.

Clearly $p M_s \subseteq M_{s'}$ and  $p M_{s'} \subseteq M_s$
by \eqref{lem-h2-i}.

 \noindent Case 1.  $s_1 \neq s_2$, and $s_1 s_2^{-1} \neq q^{\pm1}$. \\
Then $p^2$ is an invertible semisimple automorphism, and
$M = M_s \oplus p M_s$.
We have inverse equivalences of categories
$$\Rep_s S \rightleftarrows \Rep_s \Haff{2},
\qquad A \mapsto \Haff{2} \otimes_S A, \quad M_s 
\leftarrow M.
$$

\noindent Case 2.	$s_2 = q s_1$. $(q \neq 1)$ \\
In this case, $p^2$ is nilpotent, and $1 - X_1 X_2^{-1}$ is invertible.
Now, $\ker p$ is $S$-stable, and on
$\ker p$,  $T_1$ is $\frac{q-1}{1 - X_1 X_2^{-1}}$.
It follows that $\ker p$ is an $\Haff{2}$-submodule, and that
$\ker p = (\ker p)_s \oplus (\ker p)_{s'}$ is a decomposition
as $\Haff{2}$-modules.
Further, the indecomposable summands of $\ker p$ as
an $H_2^{\fin}$-module 
are the indecomposable summands as an $\Haff{2}$-module.
First suppose $q \neq -1$.
Then $T_1$ acts semisimply on $\ker p$, hence so does
$X_1 X_2^{-1}$. We suppose for simplicity that the central element
$X_1X_2$ acts semisimply.
Then $\ker p$ is a direct sum of one dimensional
$\Haff{2}$-modules.
On the subspace on which $T_1$ acts as $-1$, we have 
$X_1 = q X_2$;
on the subspace on which $T_1$ acts as $q$, we have
$X_2 = q X_1$.
Now suppose $q=-1$.
Then $(T_1 + 1)^2 = 0 = (1 - X_1 X_2^{-1})^2$,
and again there are two isomorphism classes of 
indecomposable $H_2^{\fin}$-modules: 
namely $T_1$ acts as $-1$ or as 
$\left( \begin{tabular}{rr} -1 & 1 \\  0 & -1 \end{tabular} \right)$. 

Let $A \in \Rep_{(s_1, q s_1)}$, and set
$A' \in \Rep_{( q s_1, s_1)}$ to be 
$(X_1 \leftrightarrow X_2)^* A$.
Put $A_+ = \Im(\Ind A  \to \Ind A') = p \Ind A$,
\quad $A_-= \Im(\Ind A' \to \Ind A) = p \Ind A'$,
and $A_0 = \Ind A$.
Then if $A$ is indecomposable, so are $A_0, A_+, A_-$;
and if $q= -1$ these give a complete list of 
(non-isomorphic) indecomposables.
For all $q$, the modules $(J_1 \boxtimes q J_1)_+$
and $(J_1 \boxtimes q J_1)_-$ are the distinct irreducibles.
%

\noindent Case 3.  $s = (s_1, s_1)$.\\
We have nothing general to add to the description
of these modules already given by theorem  \ref{thm-enhanced}; however
we record the following, which is proved by direct computation.
\begin{lemma}
$\Res_{\Haff{1}}^{\Haff{2}} \Ind_{\Haff{1} \otimes
\Haff{1}}^{\Haff{2}} ( q^i J_a \boxtimes q^i J_b) =
(a-1) \cdot q^i J_b +  (b-1) \cdot q^i J_a + q^i J_{a+b}$.
\end{lemma}
%

\subsection{$\Haff{3}$}
\label{sec-H3}
We list the number 
of the irreducible $\Haff{3}$-modules for which 
$\prod(s_i - q s_i) = 0$ and some of the characters of these
modules.
Note that as $\Rep_s \Haffn \isom \Rep_{s'} \Haffn$, where
$s' = (\alpha s_1, \ldots, \alpha s_n)$ for some $\alpha \in R^\times$,
we may assume $s_1 = 1$. 

\noindent $q=-1$: \\
 $\Rep_{(1,1,q)}$: 
There are 3 irreducible modules with central
character $(11q)$; two are
2 dimensional, and one is 1 dimensional. Their characters are
 $2(11q)$, $2(q11)$ and $(1q1)$.


\noindent $q^3 = 1$: \\
There are 6 irreducible representations in $\Rep_{(1,q,q^2)}$,
2 in $\Rep_{(1,1,q)}$, 2 in $\Rep_{(1,q,q)}$.
The modules in $\Rep_{(1,1,q)}$ have characters $2(11q) + (1q1)$
and $2(q11)+ (1q1)$.

\noindent $q^2 \neq 1, q^3 \neq 1$: \\
There are 4 irreducible modules in $\Rep_{(1,q,q^2)}$,
2 in $\Rep_{(1,1,q)}$, 2 in $\Rep_{(1,q,q)}$.
$\Rep_{(1,q,a)}:$
 there are 2 in $\Rep_{(1,q,a)}$,
where $a \notin \{q^{-1}, 1,q,q^2\}$;
both
$\Ind( (1q) \boxtimes a) = (a1q) + (1aq) + (1qa)$
and
$\Ind( (q1) \boxtimes a) = (aq1) + (qa1) + (q1a)$ are irreducible.
Except for $\Rep_{(1,q,q^2)}$ these modules have the same character as when
$q^3 =1$.

\subsection{$\Haff{4}$ when $q= -1$}.\\
\label{sec-H4}
$\Rep_{(1,1,1,q)}$ has three irreducibles:
$6(111q)+2(11q1) = \Ind 2(11q) \boxtimes 1$,
$2(11q1)+2(1q11) = \Ind (1q1) \boxtimes 1$,
and $6(q111) +2(1q11) = \Ind 2(q11) \boxtimes 1$.

$\Rep_{(1,1,q,q)}$ has six irreducibles, four of dimension 2 
and two of dimension 1:
$2(11qq), 2(qq11), 2(1qq1), 2(q11q), (1q1q),$ and $(q1q1)$.

\section{Bernstein-MacDonald Hopf algebra}
\label{sec-hopf}
%
Induction and restriction give $K(\Raff)$ and $K(\Rfin)$
the structure of bialgebras.
Precisely, if $M \in \Rep_q H_{a_1}^{\aff}$,
$N \in \Rep_q H_{a_2}^{\aff}$,
define multiplication
$$M \cdot N = \Ind_{H_{a_1}^{\aff} \otimes H_{a_2}^{\aff}}^{
H_{a_1 + a_2}^{\aff}} M \boxtimes N$$
and comultiplication $\Delta M = \bigoplus_{a_1 + a_2 = n}
\Delta_{a_1, a_2} M$, where $\Delta_{a_1, a_2} M =
\Res_{H_{a_1}^{\aff} \otimes H_{a_2}^{\aff}}^{H_{n}^{\aff}} M$
(and similarly for $\Rfin$).
Then
as induction and restriction are exact functors,
these descend to give functors
$$
K(\Rep_q \Haff{a_1}) \otimes K(\Rep_q \Haff{a_2})
= K(\Rep_q (\Haff{a_1} \otimes\Haff{a_2}))
\mathop{ \rightleftarrows}\limits_{\Delta_{a_1,a_2}}^{\cdot}
K(\Rep_q \Haff{a_1+a_2})
$$
%
and the properties of $\Ind$ and $\Res$ translate into the
axioms of a bialgebra; viz:
\begin{quote} 
transitivity of induction becomes associativity of multiplication;
transitivity of restriction becomes coassociativity of
comultiplication; the trivial representation of the trivial
algebra $H_0^{\aff}$ is the unit, ditto for the counit,
and the Mackey formula is the statement that $\Delta$ is an
algebra homomorphism.
\end{quote}

Similarly, define multiplication and comultiplication
on $K(\Rfin)$.
\begin{rem}
\label{rem-comult}
The comultiplication on $K(\Raff)$ is 
{\it not\/} cocommutative.
However, as we shall see,  multiplication on
$K(\Raff)$ is commutative.  
Hence the bialgebra $K(\Raff)$ is the dual to the
enveloping algebra of a Lie algebra; and Theorem \ref{uz-bi}
identifies this algebra.  
(The structure of the categories $\Rlam$ will be essential for
this identification.)
Note that $K(\Rfin)$ is both commutative and cocommutative, but
because of the previous remark, and {because of the structure of this
as an algebra over $\Z$,} it is more natural to think of this as the
{\it dual\/} to an enveloping algebra also. (That way it is a polynomial 
algebra; if we don't take the dual it is a divided power algebra.) 

\end{rem}
\begin{rem}
\label{wreath}
Structure on $\bigoplus \Rep S_n$. 
For $r > 1$, we can form embeddings of the wreath product
\begin{equation}
S_n \wr \Z_r \times S_m \wr \Z_r \hookrightarrow 
S_{n+m} \wr \Z_r
\label{wr}
\end{equation}
making $\bigoplus_{n \ge 0} \Rep (S_n \wr \Z_r)$ into a Hopf
algebra  which again classically is known
to be the $r^{th}$ tensor product of the Fock space 
$\bigoplus \Rep S_n$ (see \cite{Mac} for example).

The algebras $\Hlam{n}$ are a deformation of 
the wreath product $S_n \wr \Z_r$, and one can ask if 
the embedding \eqref{wr} 
exists for them. It is a consequence
of theorem \ref{thm-crystal} that it does {\it not}..
However, we always have an action of $K(\Rfin)$ on $K(\Rlam)$,
which is the Fock space structure when $\lambda = \Lambda_0$.
\end{rem}

This action of $K(\Rfin)$ on $K(\Rlam)$, even
when $\lambda = \Lambda_0$ is only part of the story---the coaction
of $K(\Raff)$ provides much more structure.
Eventually we will see that
the affine algebra $\sl$ acts on $K(\Rlam)$, and the 
Bernstein-MacDonald Hopf algebra structure on $K(\Rfin)$ is just the
principal realization of $\sl$. But first:

As $\evst \circ \Res = \Res \circ \evst$, we have
\begin{lemma}
\label{lemma-coalg}
The functor $K(\Rfin) \xrightarrow{\evst} K(\Raff)$ is a homomorphism of
coalgebras.
In particular, $K(\Rfin)$ is a comodule for $K(\Raff)$.
\end{lemma}
\begin{warning} $\evst$ is {\it not\/} a homomorphism of algebras.
\end{warning}

More generally, if $M \in \Rlam$, then 
$\Res_{H_{a}^{\aff} \otimes H_b^{\aff}}^{H_n^{\aff}} \evst (M)$
is an element of $\Rep_q (H_a^\lambda \otimes H_b^{\aff})$ and
\begin{lemma}
\label{lemma-comodule}
$K(\Rlam)$ is a comodule for $K(\Raff)$.
\end{lemma}
This is saying there are significantly more operations of
$K(\Rlam)$ than appear at first sight.
In particular, any exact functor $\Raff \xrightarrow{F} R$-mod gives a functor
$K(\Rlam) \to K(\Rlam)$ by composing the comodule structure 
with $\otimes_R F$
$$
K(\Rlam) \to K(\Rlam) \otimes K(\Raff) \xrightarrow{Id \otimes F}
K(\Rlam).
$$
If $\otimes_R F$ satisfies some finiteness conditions,
we obtain left and right adjoint functors $\Rlam \to \Rlam$
(which in favorable cases agree and are exact).
We will now compute some examples of this.

\section{The functors $\esti$ and $\fsti$}
\label{sec-ef}
In this section we define an action of the generators of $\sl$ on $K(\Rlam) =
\bigoplus_{n \ge 0} K(\RHlam)$.

Define functors $\est_i$ for $q^i \in \mu_q$,
$$\est_i \cc  \RqHaff \to \Rep_q H_{n-1}^{\aff}   \hspace{2ex} 
\est_i \cc  \RHlam \to \Rep H_{n-1}^{\lambda}$$
as follows.  
If $M \in \RHaff$, $\est_i M$ is the {\it generalized\/} eigenspace of
$X_n$ with eigenvalue $q^i$.
As $X_n$ commutes with $H_{n-1}^{\aff} \hookrightarrow H_n^{\aff}$,
$\est_i M$ is an $H_{n-1}^{\aff}$-module.  
Clearly $X_1$ acts in the same way on $M$ and $\est_i M$, so if we define
$\est_i M = \est_i(\evst(M))$ for $M \in \RHlam$, then we
have $\est_i M \in \Rlam$ also. 
\begin{lemma}
\label{lemma-e-exact}
The functors $\est_i \cc  \Raff \to \Raff, \est_i \cc  \Rlam \to \Rlam$
are exact.
\end{lemma}
\begin{proof}
$\est_i$ is the composite of the exact functors of  {\it restriction\/}
and the functor of  {\it generalized eigenspace\/}, which is
exact on torsion (and in particular finite dimensional)
 $R[x]$-modules.
\end{proof} 
\begin{rem}  In the abstract language favored in the previous section,
we can write
$$\est_i M = \limitm  \ker( (X_n - q^i)^m, 
\Res_{H_{n-1}^{\aff} \otimes H_{1}^{\aff}}^{H_{n}^{\aff}} M )$$
where we have identified $\Rlam H_{n}^{\aff}$ with $\RHlam$
via $\evst$.
Note that as $M$ is finite dimensional the direct limit stabilizes.
\end{rem}

As the functors $\est_i$ are exact, we may try and find their adjoints.
For  formal reasons, adjoints exist in $H_n^{\aff}$-mod and 
$H_n^\lambda$-mod, but that is no guarantee they exist in
$\RHaff$ or $\RHlam$ (i.e.~that they are finite dimensional).
Nonetheless, define functors
$\fst_i, \fhat$ from $H_{n-1}^{\aff}$-mod to $H_n^\lambda$-mod
by setting, for $M \in H_{n-1}^{\aff}$-mod,
$$ \fst_i(M) = \limitlm
\pr_\lambda ( \Ind_{H_{n-1}^{\aff} \otimes H_{1}^{\aff}}^{H_{n}^{\aff}}
( M \boxtimes q^i J_m) )$$
and
$$ \fhat(M) = {\limitm} \prhat_\lambda (
\Indhat_{H_{n-1}^{\aff} \otimes H_{1}^{\aff}}^{H_{n}^{\aff}}
( M \boxtimes q^i J_m) )$$
where $q^i J_m$ is the Jordan block of size $m$ and eigenvalue
$q^i$, $\pr$ and $\prhat$ are the left and right adjoints
to $\evst \cc  \RHlam \to \RqHaff$, 
and the direct and inverse limits are taken with respect to
the  systems
$$
q^i J_0 \hookrightarrow q^i J_1 \hookrightarrow \cdots
\hookrightarrow q^i J_m \hookrightarrow \cdots, 
\qquad
q^i J_0  \twoheadleftarrow \cdots \twoheadleftarrow q^i J_m \twoheadleftarrow \cdots
$$
given by multiplication by 
$(x-q^i)\cc  R[x] / (x-q^i)^{m-1} \to R[x] / (x-q^i)^m$, and $1 \mapsto 1$,
respectively.

Let us abbreviate $\Res_{H_{n-1}^{\aff} \otimes H_{1}^{\aff}}^{H_{n}^{\aff}}$
by $\Res$, and similarly abbreviate 
$\Ind_{H_{n-1}^{\aff} \otimes H_{1}^{\aff}}^{H_{n}^{\aff}}$ by $\Ind$.
\begin{prop}
\label{prop-stabilize}
If $N \in H_n^{\aff}$-mod, the inverse system
$$\prlam \Ind( N \boxtimes q^i J_m)$$
stabilizes after finitely many terms.
\end{prop}
\begin{proof}
If $M$ is an $H_{n+1}^{\aff}$-module generated by an $R$-subspace $W \subset M$,
then $\prlam (M)$ is an $H_{n+1}^{\lambda}$-module generated by
the image of $W$. 
In particular, if $M$ is finitely generated, then $\prlam (M)$ 
is finite dimensional (it is a quotient of the finite dimensional
vector space $ H_{n+1}^{\lambda} \otimes W$).
Now if $N \in H_n^{\aff}$-mod is generated by a subspace
$W^\prime$, $\Ind(N \boxtimes q^i J_m)$ is generated by 
$W^\prime \boxtimes q^i J_1 \isom W^\prime$
and the system $\pr(\Ind(N \boxtimes q^i J_m))$ are all quotients of
a fixed $H_{n+1}^{\lambda}$-module $H_{n+1}^{\lambda} \otimes W^\prime$,
which if $N$ is finitely generated is finite dimensional.
Hence the inverse system stabilizes in this case.

A more careful analysis shows that the system stabilizes for $m$
greater than a fixed constant which depends on $n, \lambda$ and
not on $W^\prime$.
%
This gives the proposition in general, but as we will not use it
for non-finitely generated modules,
we will omit further details.
\end{proof}
\begin{cor}
\label{cor-fstar}
If $M \in \Rep H_{n-1}^{\lambda}$, then $\fst_i M$ and $\fhat M$
are finite dimensional, i.e.~$\fst_i, \fhat$ restrict to
functors $\Rep H_{n-1}^\lambda \to \RHlam$.
\end{cor}
We have defined $\fst_i$, $\fhat$ so that
\begin{prop}
\label{prop-ef-adjoint}
The functor $\fst_i \cc  H_{n-1}^{\aff}\text{-mod} \to 
H_{n}^{\lambda}\text{-mod}$
is left adjoint to 
$\est_i \cc  H_{n}^{\lambda}\text{-mod} \to  H_{n-1}^{\aff}\text{-mod}$.
Similarly, $\fhat$ is right adjoint to $\est_i$.
\end{prop}
\begin{proof}
We prove $\fst_i$ is left adjoint; the proof for $\fhat$ is similar.
This is almost, but not quite, a formal result.
To see that, observe that if $M \in  H_{n}^{\aff}\text{-mod}$,
\begin{align*}
\est_i M &= \limitm \ker( (X_n - q^i)^m, \Res M) \\
&= \limitm \Hom_{R[X_n^{\pm 1}]}(q^i J_m, \Res M)
\end{align*}
where the second limit is taken over the system
$q^i J_0  \leftarrow q^i J_1  \leftarrow  \cdots$
used in the definition of $\fst_i$. This equals
$$ \Hom_{R[X_n^{\pm 1}]}(\limitlm q^i J_m, \Res M).  $$
if $M$ is finite dimensional, or more generally 
$R[X_n^{\pm 1}]_{(X_n -q^i)}$-torsion, but {\it not\/} in general.
(For example, if $n=1$ and $M = R[[X_n -q^i]]$ they clearly differ,
and indeed $\est_i$ is {\it not\/} exact on
the category of all $H_n^{\aff}$-modules.)

However, if $M \in H_n^\lambda$-mod, then as $H_n^\lambda$ is finite 
dimensional, these are equal, and the direct and inverse limits above
stabilize after finitely many terms (for $m \ge \dim H_n^\lambda$, to
be crude).
Hence if $M \in H_n^\lambda$-mod, and $N \in H_{n-1}^{\aff}$-mod, then
%
\begin{align*}
\Hom_{ H_{n-1}^{\aff}}(N, \evst(\est_i M)) &=
\Hom_{ H_{n-1}^{\aff}}(N, \limitm  \Hom_{R[X_n]}(q^i J_m, \Res \evst M ))\\
&= \limitm \Hom_{ H_{n-1}^{\aff}}(N,\Hom_{R[X_n]}(q^i J_m, \Res \evst M ))\\
\intertext{as the direct limit stabilizes after finitely many terms, and this equals} 
&= \limitm \Hom_{ H_{n-1}^{\aff} \otimes H_1^{\aff}}
(N \boxtimes q^i J_m, \Res \evst M )\\
&= \limitm \Hom_{ H_{n}^{\aff}} (\Ind(N \boxtimes q^i J_m), \evst M ),\\
\intertext{as $\Ind$ is left adjoint to $\Res$.  As $\pr_\lambda = \pr$ is
left adjoint to $\evst$, we get}
&= \limitm \Hom_{ \Hlam{n}} (\prlam \Ind(N \boxtimes q^i J_m),M) \\
&= \Hom_{ \Hlam{n}}(\limitlm \prlam \Ind( N \boxtimes q^i J_m),M) \\
&= \Hom_{H_n^\lambda}( \fst_i N, M)
\end{align*}
%
where we have used, once again, the fact that this limit of
Hom's stabilizes after finitely many terms.
\end{proof}

In particular, restricting to finite dimensional modules $\Rlam$, we get
\begin{cor}
\label{cor-rtlftadjoint}
$\fst_i \cc  \Rep H_{n-1}^\lambda \to \RHlam$ is left
adjoint to $\est_i \cc  \RHlam \to \Rep H_{n-1}^\lambda$, and
$\fhat \cc  \Rep H_{n-1}^\lambda \to \RHlam$ is right adjoint.
\end{cor}
Note that the proof above shows for any $N \in  H_{n-1}^{\aff}$-mod,
and any $M\in \RHlam$, that
$$\Hom(\fst_i N, M) = \Hom( \prlam \Ind(N \boxtimes \limitlm q^i J_m), M)
$$
so by uniqueness of adjoints 
$$
\fst_i N = \prlam \Ind(N \boxtimes \limitlm q^i J_m).
$$
(But it is certainly {\it not\/} true that
$\fst_i N = \prlam \Ind(N \boxtimes \limitm q^i J_m)$
for arbitrary $N \in H_{n-1}^{\aff}$-mod, or even for arbitrary $N
\in \Rep H_{n-1}^\lambda$. Obviously.)

Recall that $\RHaff = \bigoplus \Rep_s H_n^{\aff}$,
where the sum is over $(R^\times)^n /S_n$, the orbits of
the symmetric group on $(R^\times)^n$.
A module $M$ is in $ \Rep_s H_n^{\aff}$ if $M$, when considered as
a module for $Z(H_n^{\aff})$, has support a single orbit of
$S_n$ on $(R^\times)^n$.
Further, if $M$ is indecomposable and 
$s = (s_1, \ldots, s_n) \in (R^\times)^n$ is such that
the weight space $M_s = \{ m \in M \mid X_i \cdot m = s_i m, 1 \le i \le n \}$
is non-zero, then $M \in \Rep_s  H_n^{\aff}$.
(Conversely, if $M$ is indecomposable, there is always some
$s \in (R^\times)^n$ such that $M_s \neq 0$.)

If $s = (s_1, \ldots, s_n)$, and $q^i \in \{s_1, \ldots, s_n\},$
say $s_n = q^i$, let's write $s \setminus  q^i$ for the unique
orbit in $(R^\times)^{n-1}/S_{n-1}$ obtained by deleting $q^i$ from
the list.  Then lemma \ref{lemma-e-exact} has a refinement (which needs
no proof).
\begin{lemma}
\label{lemma-s-q}
$\est_i$ is a functor $\Rep_s H_n^{\aff} \to \Rep_{s \setminus q^i}
 H_{n-1}^{\aff}$,
and $\est_i M = 0$ if $q^i \not\in  \{s_1, \ldots, s_n\}.$
\end{lemma}
Dually, write $s + q^i$ for the orbit of $S_{n+1} \cdot
(s_1, \ldots, s_n, q^i) \in (R^\times)^{n+1}/S_{n+1}.$
We can now refine the previous proposition to:
\begin{prop}
\label{prop-adjef}
$\fst_i \cc  \Rep_s H_{n-1}^\lambda \to \Rep_{s+q^i} H_{n}^\lambda$
is left adjoint to 
$\est_i \cc  \Rep_{s+q^i} H_{n}^\lambda \to \Rep_s H_{n-1}^\lambda$,
and $\fhat$ is right adjoint.
\end{prop}
\begin{proof}
Let $M$ be an indecomposable $H_{n-1}^\lambda$-module, and $v$ a non-zero
element in the weight space $M_s$ for some $s \in (R^\times)^n$.
Then all vectors in $Rv \boxtimes q^i J_m$ are in the generalized 
eigenspace with eigenvalue $s +q^i$, and hence all subquotients 
of $\Ind(\evst M \boxtimes q^i J_m)$ are in $\Rep_{s+q^i} H_{n}^{\aff}$.
But if $N \in \Rep_s H_{n}^{\aff}$, then $\pr_\lambda N \in \Rep_s H_n^\lambda$,
and so $\fst_i M \in \Rep_{s+q^i} H_{n}^\lambda$. 
(The proof for $\fhat$ is similar.) 
\end{proof}
\begin{lemma}
If $M \in \RHlam$, then
\begin{enumerate}
\item[(i)]  $\Res_{H_{n-1}^\lambda}^{H_{n}^\lambda} M \isom
\bigoplus_{i \in \mu_q} \est_i M$    \label{e=res} \\
\item[(ii)]  $\Ind_{H_{n}^\lambda}^{H_{n+1}^\lambda} M \isom
\bigoplus_{i \in \mu_q} \fst_i M \isom \bigoplus_{i \in \mu_q} \fhat M$ 
\label{f=ind}
\end{enumerate}
\end{lemma}
\begin{proof}
For $M \in \RHlam$, $X_n$ acts on $M$ with eigenvalues in
$\mu_q$, so \eqref{e=res} is immediate from the definition of
$\est_i$.  But then \eqref{f=ind} follows as both 
$\sum \fst_i$ and $\Ind_{H_{n-1}^\lambda}^{H_{n}^\lambda}$
are left adjoint to $\Res_{H_{n-1}^\lambda}^{H_{n}^\lambda} \cc 
\RHlam \to \Rep H_{n-1}^\lambda$, and left adjoints are 
unique when they exist.
\end{proof}
\begin{cor}
\label{cor-f-exact}
As functors $\Rlam \to \Rlam$, the functors $\fst_i, \fhat$
satisfy
\begin{enumerate}
\item[(i)]  
$\fst_i \isom \fhat$ \\
\item[(ii)] 
$\fst_i$ is exact.
\end{enumerate}
\end{cor}
\begin{proof}
As $\fst_i, \fhat$ are direct summands of the
exact functor $\Ind_{H_{n}^\lambda}^{H_{n+1}^\lambda}$, 
they are exact.
Proposition \ref{prop-adjef} identifies them as the { same\/}
direct summands.
\end{proof}

\begin{rem} We could also use this as a definition of $\fst_i$, but then
we would not see the uniform dependence on $\lambda$, or be able 
to prove the theorems of the next section.
\end{rem}

\begin{rem}
Note that for
$N \in \Rep H_{n-1}^\lambda $,
$$
\Ind (\evst N \boxtimes q^i J_m) \neq \Indhat (\evst N \boxtimes q^i J_m).
$$
\end{rem}

Because $\fst_i$ is left 
adjoint to an exact functor, or because
restriction and induction take free modules to free modules, we observe
\begin{lemma}
If $M \in \Rlam$ is projective, then so are $\est_i M$ and $\fst_i M$.
\end{lemma}
%
\subsection{Divided Powers}
\label{sec-divided}
We can do slightly better.
Fix $i \in \mu_q$.  
For each $n \ge 1$ we define the {\it divided powers\/}
of $\esti$ and $\fsti$ as exact functors $\Rlam \to \Rlam$.

First, define $\edst\cc  \Rep \Haff{a+n} \to \Rep \Haff{a}$
as follows
$$
\edst (M) = \mathop{\varinjlim}\limits_P \: \Hom_{\Haff{a} \otimes \Haff{n}}
(\Haff{a} \boxtimes P, \Res_{\Haff{a} \otimes \Haff{n}}^{\Haff{a+n}} M)
$$
where the limit is taken over the small category whose objects are the
finite dimensional approximations $P$ to a projective module for $\Haff{n}$
which surjects onto $\Ind_{\Haff{1} \otimes \cdots \otimes \Haff{1}}^{\Haff{n}}
(q^i J_1 \boxtimes \cdots \boxtimes q^i J_1)$.
More precisely, the objects of this category are the diagrams
$$ {\mathcal P} \surj P \surj \Ind_{\Haff{1} \otimes \cdots \otimes \Haff{1}}^{\Haff{n}}
(q^i J_1 \boxtimes \cdots \boxtimes q^i J_1)$$
where $P$ is a finite dimensional $\Haff{n}$-module and ${\mathcal P}$ is 
a projective module for $\Haff{n}$.
Morphisms from $P$ to $P^\prime$ are commutative diagrams such that
all maps are surjective.

Then define $\edst \cc  \Rep H_{a+n}^\lambda \to \Rep H_a^\lambda$
as $\edst \circ \evst_\lambda$, and define
$\fdst \cc  \Rep_q \Haff{a} \to \Rep H_{a+n}^\lambda$ by
$$
\fdst (M) =  \mathop{\varprojlim}\limits_P \pr_\lambda 
\Ind_{\Haff{a} \otimes \Haff{n}}^{\Haff{a+n}} (M \boxtimes P).
$$
We omit the proof that these are reasonable definitions---that
they are exact functors, that $\fdst$ is both left and
right adjoint to $\edst$, that they preserve the properties of
being finite dimensional, etc.  

More generally given an irreducible module
$L \in \RHaff$ we get functors
$$
\Delta_L\cc  \Rep \Haff{a+n} \to \Rep \Haff{a}, \qquad
m_L \cc  \Rep \Haff{a} \to \Rep H_{a+n}^\lambda
$$
by mimicking the above construction---just replace
$\Ind(q^i J_1 \boxtimes \cdots \boxtimes q^i J_1)$
above with $L$.
If $L$ is the module $q^i J_1$ for $\Haff{1}$,
$\Delta_L = \esti = {\e_i^{(1)*}}$ and 
$m_L = \fsti = {\f_i^{(1)*}}$.
In general (and even for $L = q^i J_1$), these functors
only have good properties in the Grothendieck group.
It is clear $\Delta_L$ is the composite of the
comodule action of $K( \Raff)$ on $K(\Rlam)$ with the
function $\delta_L$ of 
``multiplicity of $L$ in the Jordan-Holder series''
$$
\Delta_L\cc  K(\Rlam) \xrightarrow{\Delta} 
K(\Rlam) \otimes K(\Raff) \xrightarrow{1 \otimes [L:\quad ]}
K(\Rlam) \otimes_\Z \Z = K(\Rlam)
$$
(see section \ref{sec-grothgroup} 
for the notation).
In order to describe $m_L$ on the level of the
Grothendieck group we must introduce the
Shapovalov inner product in $K(\Rlam)$; 
then $m_L$ will be the adjoint of $\Delta_L$ with respect
to this inner product. This is done in section \ref{sec-shap}. 
\section{crystal graph}
\label{sec-crystal}

In this section we will study certain ``highest order'' approximations
to $\esti$ and $\fsti$.  These approximations make sense for both
cyclotomic and affine Hecke algebras, unlike $\fsti$ itself, which is
not defined for $\Raff$. We summarise the properties of these
operators in a combinatorial structure, the crystal graph.

There are three main results in this section. 
The first  is theorem \ref{thm-multone}, which is a strong
multiplicity one theorem for restriction and induction.
This tells us that the crystal graph is well defined.

The other two results are theorems \ref{thm-epsilon} and theorem
\ref{thm-nu}.  These describe the close connection between the the
crystal graph structure and the representation theory of the Hecke
algebra.

Theorem \ref{thm-epsilon} shows that the integers $\epsi(M)$ have
various interpretations.  By definition, if $M$ is irreducible
$\epsi(M)$ is the length of the longest chain of $q^i$'s that end the
spectrum of $M$. But we show that it is also the maximum size of a
Jordan block of $X_n$ with eigenvalue $q^i$ on the module $M$, i.e.~
that it measures the failure of semsimplicity of the action of $X_n$.
It is also the dimension of $\Hom(\esti M,\esti M)$---another subtle
measure of lack of semisimplicity, as well as the multiplicity of the
cosocle of $\esti M$ in $\esti M$. These last two statements show that the
cosocle of $\esti(M)$ fits into $\esti(M)$ in a uniserial chain of
length $\epsi(M)$ in as simple a way as possible.

The analogous results for $\fsti$ is theorem \ref{thm-nu}. This is
much harder, and the particular proof we give requires the results of
the next few sections.  The asymmetry is a nice shadow of the fact
that whereas the dual of a finite dimensional lowest weight
${\mathfrak sl}_2$ module is again a lowest weight module, the dual of
an integrable lowest weight $\sl$-module is a {\it highest} weight
module. Nonetheless, the ${\mathfrak sl}_2$ structure is the relevant
one, and we eventually show $\nui(M)$ may also be read off the
spectrum of $M$.

The results of section \ref{sec-shap} show that $\epsi$, $\nui$  also
admit an interpretation in terms of the structure of projective modules.

\subsection{First properties}
To save repetition in notation, we will allow  $\Hlam{n}$ to denote the affine
Hecke algebra as well as the cyclotomic Hecke algebra, i.e.~we allow
$\lambda$ to be either the symbol ``$\aff$'', or a function
$\lambda \cc  \mu_q \to \Z_+$ with $\sum \lambda_i < \infty$.
We say $\lambda$ is {\it affine\/} or cyclotomic, as appropriate.
For consistency in notation, define $\pr_{\aff}$, $\evst_{\aff}$ to be
the identity functors from $\RqHaff \to \RqHaff$.

With this understood we write
$\B$ for the disjoint union of the set of isomorphism
classes of irreducible representations on $H_n^\lambda$, for $n \ge 0$.
Let $\Z\B$ be the free abelian group on the set $\B$, and define
a non-degenerate symmetric pairing on $\Z\B$ by making the elements of $\B$
an orthonormal basis.
We will identify $\Z_+\B$ with isomorphism classes of semisimple
$H_n^\lambda$-modules, and the symmetric pairing with
$\dim_R \Hom(\quad , \quad )$.
(There is a canonical isomorphism of abelian groups---without the
symmetric pairing---from the Grothendieck group $K(\Rlam)$ to
 $\Z\B$, but this will not be so useful for what follows.)

Now, define for $M \in \RHlam$
\begin{gather*} \ft(M) = \pr_\lambda \cosoc 
\Ind_{\Haff{n} \otimes \Haff{1}}^{\Haff{n+1}} 
(\evst_\lambda M \boxtimes q^i J_1) \\
\et(M) = \soc \Hom_{R[X_n^{\pm 1}]}(q^i J_1, \Res_{\Haff{n-1} \otimes
 \Haff{1}}^{\Haff{n}} \evst_\lambda M).
\end{gather*}
As $\et(M)$ is semisimple and $\et(M \oplus M^{\prime}) = \et(M)
\oplus \et(M^{\prime})$, $\et$ defines a $\Z_+$-linear
operator on isomorphism classes of semisimple $\Hlam{n}$-modules, and hence
an operator $\et \cc  \Z \B \to \Z \B$.
Similarly, $\ft$ defines an operator $\Z \B \to \Z \B$, as 
$\pr_\lambda$ of a semisimple module is still semisimple.
These are adjoint operators with respect to the
inner product on $\Z \B$, i.e.:
\begin{lemma}
If $N$ is a semisimple $\Hlam{n-1}$-module, and $M$ a semisimple $\Hlam{n}$-module,
$$
\Hom_{\Hlam{n-1}}(N, \et M) = \Hom_{\Hlam{n}}(\ft N, M).
$$
\end{lemma}
\begin{proof}
As $N$ is semisimple,

\begin{gather*}
\Hom_{\Haff{n-1}}(N, \soc \Hom_{R[X_n^{\pm1}]}(q^i J_1, \Res \evst M))
\qquad { } \qquad {  } \qquad \quad\\
\qquad \quad
{\begin{align*}
&= \Hom_{\Haff{n-1}}(N,\Hom_{R[X_n^{\pm1}]}(q^i J_1, \Res \evst M)) \\
&= \Hom_{\Haff{n-1} \otimes \Haff{1}}(N \boxtimes q^i J_1, \Res \evst M)
\\
&=  \Hom_{\Haff{n}}(\Ind(N \boxtimes q^i J_1),\evst M ) \\
\intertext{by adjointness. Now, as $\evst M$ is semisimple, this is}
&=  \Hom_{\Haff{n}}(\cosoc\Ind(N \boxtimes q^i J_1), \evst M )\\
&=  \Hom_{\Hlam{n}}(\pr_\lambda \cosoc \Ind(N \boxtimes q^i J_1), M).
\end{align*}}\\
\end{gather*}
\end{proof} 
As with the operators $\fsti$ the order
of the operations defining $\ft$ is not crucial.
We have:
\begin{lemma}
$\ft(M) = \cosoc (\pr_\lambda \Ind_{\Haff{n-1} \otimes \Haff{1}}^{\Haff{n}}
(\evst_\lambda M \boxtimes q^i J_1))$.
\end{lemma}
\begin{proof}
We may finish the above proof differently, noticing $\evst M$
is still semisimple, to get 
\begin{multline*}
\Hom_{\Hlam{n-1}}(N, \et M) = 
\Hom_{\Hlam{n}} (\pr_\lambda \Ind(N \boxtimes q^i J_1), M) \\
= \Hom_{\Hlam{n}}(\cosoc \pr_\lambda \Ind(N \boxtimes q^i J_1), M)
\end{multline*}
It follows that both $\ft$
and this new operator define adjoints to 
$\et \cc  \Z \B \to \Z \B$.
But the symmetric pairing $\Hom ( \;, \; )$ is non-degenerate,
and so these two operators must be equal.
\end{proof}
In fact, we will show in theorem \ref{thm-ft-cosocf*} that we 
may even replace $J_1$ with $J_m$ in the definitions of
$\et$ and $\ft$. 

We have implicitly used the map
$$
\evst\cc  \Z \B \hookrightarrow \Z B_\aff 
$$
induced by $\evst\cc  \RHaff \to \RHlam$.
The definitions of $\esti $ and $\fsti$ make the following lemma obvious.
\begin{lemma}
\label{lemma-evei-eiev}

\begin{enumerate}
\item[(i)] $\et \evst M = \evst \et M$, for $M \in \B$. \\
\item[(ii)] If $M \in \B$, $\evst \ft M$ is a direct summand of $\ft(\evst M)$.
\end{enumerate}
\end{lemma}
\begin{thm}
\label{thm-multone}
If M is an irreducible $\Hlam{n}$-module, then $\et M$ and 
$\ft M$ are either irreducible or zero.
Further, if $N \neq 0$, then $\et M = N$ if and only if
$M = \ft N$.
\end{thm}
The theorem will be proved after proposition \ref{prop-epsilon} below.

As a result, we can summarize the operators $\et$
and $\ft$ in the datum of an oriented graph, with
edges labelled by the elements of $\mu_q$.
The vertices of the graph are the elements of $\B$,
and $M$ is joined to $N$ by an arrow colored by 
$i \in \mu_q$
$$
M \xrightarrow{i} N
$$
if $\et M = N$; equivalently if $N = \ft M$.
This datum $(\B, \et, \ft)$ is called a {\it crystal
graph\/}, after Kashiwara. Note that as a consequence of
the lemma \ref{lemma-evei-eiev}, $\evst$ induces an inclusion of crystal
graphs, $\B \to B_\aff$.
We record one immediate formal consequence of this fact:
\begin{cor}
If $M \in B_\aff$, $(\prlam \ft)^k M = \prlam(\ft^k M)$.
\end{cor}
\begin{proof}
By induction on $k$, $(\pr \ft)^k M  = \pr \ft \pr \ft^{k-1}M$,
and $\pr \ft^{k-1}M$ is either $0$ or 
$\ft^{k-1}M$, so $(\pr \ft)^k M$ is either $0$ or $\pr \ft^k M$.
So if $\pr \ft^k M =0$, then
$(\pr \ft)^k M = 0$ also, and if $\pr \ft^k M \neq 0$, it equals
$\ft^k M$ (by the theorem), so
$\et(\pr \ft^k M) = \ft^{k-1}M \in \RHlam$,
and so $\ft^{k-1}M = \pr \ft^{k-1}M$.
Hence $(\pr \ft)^k M = \pr \ft^k M$ here also.
\end{proof}

Define, for $M \in \B$
$$
\epsi(M) = \max \{ n \ge 0 \mid  \et^nM \neq 0 \} $$
$$ \nui(M)    = \max\{ n \ge 0 \mid  \ft^nM \neq 0 \} $$
so that $\epsi (M) \in \Z_+$ (and indeed if 
$M \in \RHlam$, $\epsi (M) \le n$), and 
$\nui (M) \in \Z_+ \cup \{ \infty \}$.
Note that if $\lambda = \aff$ then 
$\nui (M) = \infty$ always. We will see in theorem \ref{thm-nu}
that if $\lambda$ is
cyclotomic then $\nui(M)$ is always finite.

\subsection{ Detailed study of the crystal graph}
\label{sec-detailedcrystal}
We now start to seriously study the crystal graph.
In order to lighten notation, let's agree to write
$\Ind = \Ind_{a_1, \ldots, a_k}^n = 
\Ind_{\Haff{a_1} \otimes \Haff{a_2} \otimes \cdots \otimes \Haff{a_k}}^{\Haff{n}}$
for the induction functor between modules for affine Hecke
algebras, and to omit indices when this causes no confusion;
and similarly for restriction.

Let us also write
$$
q^i K_n = \Ind_{\Haff{1} \otimes \cdots \otimes \Haff{1}}^{\Haff{n}}(q^i J_1 \boxtimes
\cdots \boxtimes q^i J_1)
$$
for the unique irreducible $\Haff{n}$-module with central character
$(q^i \cdots q^i)$, and, if 
$M \in \Rep \Haff{a}$ and $N \in \Rep \Haff{n+a}$
$$\Hom_{\Haff{a}}(M,N)$$
for the $\Haff{n}$-module     $\Hom_{\Haff{n} \otimes \Haff{a}}
(\Haff{n} \boxtimes M, \Res_{\Haff{n} \otimes \Haff{a}}^{\Haff{n+a}} N)$.
So in our lighter notation,
$\et M = \soc
\Hom_{\Haff{1}}(q^i J_1, M)$.
\begin{lemma}
Let $M \in \Rep \Hlam{n+a}$. The following are equivalent
$$
(i)\quad \epsi(M) \ge n \qquad  (ii) \quad \et^n M \neq 0 \qquad
(iii) \quad \Hom_{\Haff{n}}(q^i K_n , M) \neq 0.
$$
\end{lemma}
\begin{proof}
Clearly, (i) and (ii) are equivalent to
the existence of a nonzero map of $R [ X_{a+1}^{\pm 1}, \ldots , X_{a+n}^{\pm 1}]$-modules
$$
q^i J_1 \boxtimes \cdots \boxtimes q^i J_1 \to \Res_{\Haff{a} \otimes 
\Haff{1} \otimes \cdots \otimes \Haff{1}}^{\Haff{a} \otimes \Haff{n}} M.
$$
By Frobenius reciprocity, this is equivalent to the existence
of a nonzero map of $\Haff{n}$-modules
$$
q^i K_n \to \Res_{\Haff{a} \otimes \Haff{n}}^{\Haff{a+n}} M.
$$
\end{proof}
The  crucial technical result we will need to prove the theorem is:
\begin{prop}
\label{techprop}
Let $M \in \Rep \Haff{a}$, and suppose $\epsi(M) = 0$.
Then the exact sequence of $\Haff{a} \otimes \Haff{n}$-modules
$$
0 \to M \boxtimes q^i K_n \to \Res_{a, n}^{a+n} \Ind_{a, n}^{a+n} (M \boxtimes q^i K_n) 
\to \A \to 0
$$
splits.  Moreover, for every subquotient $A$ of $\A$,
$$
\Ext_{\Haff{a} \otimes \Haff{n}}
(A, \Haff{a} \boxtimes q^i K_n) = 0.
$$
\end{prop}
\begin{proof}
By the Mackey formula, $\A$ admits a filtration whose graded
pieces are 
$$
\Gamma_w = \Ind_{(a,n) \bigcap w(a,n)}^{(a,n)} w^{-1} \Res_{(a,n) \cap w(a,n)}^{(a,n)}
(M \boxtimes q^i K_n)
$$
where $w$ runs over representatives for all the
 cosets $S_{(a,n)} \left\backslash S_{a+n} \right/ S_{(a,n)}$ 
{\it except\/} for the coset $S_{(a,n)}$.

Now consider $\Gamma_w$ as a module for $R [ X_{a+1}^{\pm 1}, \ldots , X_{a+n}^{\pm 1}]$,
and suppose $s = (s_{a+1}, \ldots,\linebreak[3] s_{a+n}) \in (R^\times)^n$
is in its support; i.e.~$\Gamma_w [s] \neq 0$.
Fix $w \in S_{(a,n)} \left\backslash  S_{a+n} \right/ S_{(a,n)}$, not
the identity double coset $S_{(a,n)}$.
By the {\it shuffle lemma\/} \ref{shufflelemma}, there must
be some $s_\gamma$, $a < \gamma \le a+n$, 
such that the $s_\gamma$-weight space of $X_a$ on $M$ is
nonzero, i.e.~%
$\{ m \in M \mid (X_a - s_\gamma)^{\dim M} \cdot m = 0 \} \neq 0$.
As $\et M = 0$, $s_\gamma \neq q^i$.
So we have shown there is some $\gamma$ , $a < \gamma \le a+n$,
with $s_\gamma \neq q^i$.

It follows that $\Gamma_w$, when considered as a module for
$Z(\Haffn)$ has support disjoint from 
$(q^i \cdots q^i)$.
As this is the support of $q^i K_n$, we see that all
subquotients of $\Gamma_w$, considered as an
$\Haffn$-module, are in
different blocks from $q^i K_n$. The proposition follows.
\end{proof}
Write $\ftd M = \cosoc \Ind (M \boxtimes q^i K_n)$,
for $M \in \Rep \Haff{a}$.
\begin{prop}
\label{prop-epsilon}
Suppose $M \in \Rep \Haff{a}$ is irreducible, and $\epsi (M) = 0$. Then
\begin{enumerate}
\item[(i)]
\label{pone} 
$\ftd M$ is irreducible.
\item[(ii)] 
\label{ptwo}
If $A$ is an irreducible subquotient of the kernel of the map
$$
\Ind(M \boxtimes q^i K_n ) \twoheadrightarrow \cosoc \Ind(M \boxtimes q^i K_n )
$$
then $\epsi (A) < n$ Further, $\epsi(\ftd M) =n$.
\item[(iii)] 
\label{pthree}
If $M,N \in \Rep \Haff{a}$ are irreducible, and $\epsi (M) = 0$
then $\ftd M = \ftd N$ implies $M = N$.
\item[(iv)] 
\label{pfour}
$\ftd \ft^{(m)} M = \ft^{(n+m)} M$.  In particular, $\ftd M = \ft^n M$.
\end{enumerate}
\end{prop}
\begin{proof}
(i) 
As $\ftd M$ is semisimple, we will know it is irreducible
if we show $d = \dim \Hom (\ftd M, \ftd M)$ equals 1.
Note that $\cosoc (X) = 0 $ only if $X = 0$, so 
$\ftd M$ is nonzero, and we must show $d \le 1$.

Write $S$ for the kernel
\begin{gather}
\label{exactseq}
0 \to S \to \Ind(M \boxtimes q^i K_n) \twoheadrightarrow \ftd M
= \cosoc \Ind(M \boxtimes q^i K_n ) \to 0. 
\end{gather}
As $\Res = Res_{n,a}^{n+a}$ is exact, we have an exact sequence
$$
0 \to \Res S \to \Res \Ind(M \boxtimes q^i K_n)  \twoheadrightarrow 
\Res \ftd M \to 0,
$$
and by the technical proposition \ref{techprop}
$\Res \Ind(M \boxtimes q^i K_n) = M \boxtimes q^i K_n \oplus \A$.
So $\Res \ftd M = \overline{M \boxtimes q^i K_n} \oplus \overline{\A}$,
where  $\overline{M \boxtimes q^i K_n}$ is a quotient of $M \boxtimes q^i K_n$,
and $\overline{\A}$ is a quotient of $\A$. Now,
\begin{align*}
\Hom(\ftd M, \ftd M) &= \Hom (\Ind(M \boxtimes q^i K_n), \ftd M)\\
&= \Hom(M \boxtimes q^i K_n, \Res \ftd M) \\
\intertext{as $\ftd M$ is semisimple, and $\Ind$ adjoint to $\Res$.
 But this equals}
&= \Hom(M \boxtimes q^i K_n,\overline{M \boxtimes q^i K_n})
\end{align*}
%
by the technical proposition \ref{techprop}.
As $M \boxtimes q^i K_n$ is irreducible, either 
$\overline{M \boxtimes q^i K_n}$ is zero, and $d=0$,
or $\overline{M \boxtimes q^i K_n} = M \boxtimes q^i K_n$,
and $d=1$. We have already observed $d \neq 0$,
so $\ftd M $ is irreducible, nonzero.

(ii) 
As the sequence \eqref{exactseq} is exact, this also shows
that $\Res S$ embeds into $\Res \A$.
Now let $A$ be an irreducible subquotient of $S$.
Then to show $\epsi (A) < n$, it suffices to
show $\Hom_{\Haffn}(q^i K_n, \Res_{a,n}^{a+n} A)$ is zero.
But $\Res A$ is a subquotient of $\Res \A$, and
the technical proposition \ref{techprop}  gives the result.
As $\Res \ftd M = M \boxtimes q^i K_n \oplus \overline{\A}$ it is clear 
that  $\epsi(\ftd M) =n$.

(iii) 
Suppose $N \in \Rep \Haff{a}$ is irreducible.
Then by semisimplicity of $\ftd M$ and adjunction, we get
\begin{eqnarray*}
\Hom (\ftd N, \ftd M) &=& \Hom( N \boxtimes q^i K_n, \Res \ftd M)\\
&=& \Hom( N \boxtimes q^i K_n, M\boxtimes q^i K_n), \text{ as before}\\
&=& \Hom(N, M), \text{ as $N, M$ semisimple.}
\end{eqnarray*}
So $\ftd N = \ftd M \iff N = M$.

(iv) 
To show  (iv), 
observe that as we have a surjection
$$
\Ind( M \boxtimes q^i K_m) \boxtimes q^i K_n
\twoheadrightarrow \cosoc \Ind( M \boxtimes q^i K_m)\boxtimes q^i K_n
= \ft^{(m)} M \boxtimes q^i K_n
$$
we also get a surjection
$$
\Ind( M \boxtimes q^i K_{m+n}) = \Ind  ( \Ind( M \boxtimes q^i K_m) \boxtimes q^i K_n)
\twoheadrightarrow \Ind( \ft^{(m)}M \boxtimes q^i K_n).
$$
Given a surjective map $X \twoheadrightarrow Y$,
we get an induced surjection on cosocles, $\cosoc X \twoheadrightarrow
\cosoc Y$, so a surjective map
$$
\ft^{(m+n)}M \twoheadrightarrow \ftd(\ft^{(m)}M).
$$
As $\ft^{(m+n)}M$ is irreducible, this is an isomorphism.
\end{proof}
We now prove the theorem \ref{thm-multone}. 
\begin{proof}
It is clearly enough to prove it for $\Haffn$-modules, as lemma
\ref{lemma-evei-eiev} implies the result for cyclotomic Hecke algebras.

Suppose $\overline{M} \in B_\aff$ is an irreducible $\Haff{m}$-module,
and write $n= \epsi(\overline{M})$.
Then if $M$ is semisimple,
$$
\Hom(M, \et^n \overline{M}) = \Hom(\ft^n M, \overline{M}).
$$
So if $M$ is a simple summand of $\et^n \overline{M}$, then
$\et M = 0$ and the proposition \ref{prop-epsilon} shows 
$\ft^n M = \overline{M}$.
If $M^{\prime}$ is another simple summand of $\et^n \overline{M}$,
then $\ft^n M^{\prime} = \ft^n M$, so $M^{\prime} = M$.
It follows that $M = \et^n \overline{M}$.

We have just seen that {\it every\/} irreducible module
is of the form $\ft^n M$, where $\epsi(M)=0$.
It follows that $\ft N$ is irreducible, for all irreducible
modules $N$.
Finally, let $X$ be a summand of $\et \ft^n M$, where
$\epsi(M) = 0$, $M$ is irreducible.
Then $X = \ft^m N$, for some $N$ with $\epsi(N) = 0$,
as $\Hom (X, \et \ft^n M) = \Hom (\ft X, \ft^n M)$.
It follows that $\ft X = \ft^n M$, so $\ft^{m+1} N = \ft^n M$,
whence $m = n-1$ and $N=M$. So $\et \ft^n M = X = \ft^{n-1}M$,
as desired.  
\end{proof}
The previous theorem used the result that $q^i K_n$ is an irreducible
in $\Rep_s \Haffn$, where $s=(q^i \cdots q^i )$.
The next theorem uses the fact that this is the {\it only\/}
irreducible in this block.
\begin{thm}

\label{thm-ft-cosocf*}
\begin{enumerate}
\item[(i)] For any irreducible $M \in \RHlam$, and $m >0$
$$
\ft(M) = \pr_\lambda \cosoc \Ind_{\Haffn \otimes \Haff{1}}^{\Haff{n+1}}
(\evst_\lambda M \boxtimes q^i J_m).
$$
In particular, if $\lambda$ is cyclotomic, then $\ft(M) = \cosoc \fsti(M)$.
\item[(ii)]  For any irreducible $ M \in \RHlam$, and $m >0$
$$
\et (M) = \soc \Hom_{R[X_n^{\pm 1}]} (q^i J_m, \Res_{\Haffn \otimes \Haff{1}}^{\Haff{n+1}} \evst_\lambda M).
$$
In particular, $\et(M) = \soc \esti (M)$.
\end{enumerate}
\end{thm}
\begin{proof}
It is enough to prove (i).
For, if $M$ and $N$ are semisimple $\Hlam{n}$-modules, then by
(i) and adjointness 
$\Hom(\ft M, N) = \Hom(\cosoc \fsti M , N) = \Hom (\fsti M, N) = \Hom(M, \esti N) = \Hom(M,\soc \esti N)$
and so both $\et$ and $\soc \esti$ are adjoint to 
$\ft\cc  \Z \B \to \Z \B$.
As the symmetric pairing $\Hom(\quad,\quad)$ is non-degenerate,
they must be equal.

So we prove (i).  In fact, we prove a stronger statement.
Let $P \in \RqHaff$ be any indecomposable module with cosocle
 $q^i K_n$, so $P \twoheadrightarrow q^i K_n$.
Define $\f_P(M) = \cosoc \Ind(M\boxtimes P)$, for any $M \in \RqHaff$.
Clearly $\f_P(M) \twoheadrightarrow \ftd(M)$.
We will show that if $\epsi (M) = 0 $, $\f_P(M) =\ftd(M) = \ft^n(M)$
is irreducible.
Now, the socle filtration of $P$ has subquotients isomorphic
to $q^i K_n$.  It follows that we can filter 
$\Ind(M \boxtimes P)$ with subquotients isomorphic to
$\Ind(M \boxtimes q^i K_n)$. Hence $\f_P(M) = \cosoc \Ind(M \boxtimes P)$
is some summand of some copies of $\cosoc \Ind(M \boxtimes q^i K_n) 
= \ftd M$ (for any filtered module $Z = \cup Z_j$, $\cosoc(Z)$ is
a summand of $\oplus \cosoc(Z_j/Z_{j-1})$).
But 
\begin{multline*}
\Hom(\f_P M, \ftd M) = \Hom (\Ind(M \boxtimes P), \ftd M) 
= \\ \Hom(M \boxtimes P, \Res \ftd M)
= \Hom(M \boxtimes P, M \boxtimes q^i K_n) = R.
\end{multline*}
It follows that $\f_P M = \ftd M$.

Note that this implies if $\cosoc(P) = q^i K_n$, and 
$\cosoc(P')  = q^i K_{n^{\prime}}$, then
$\f_P \f_P'  = \ft^n \ft^{n^{\prime}} = \ft^{n+{n^{\prime}}} = 
\f_{\Ind(P \boxtimes P')}$.

Now, if $\overline{M}$ is any simple $\Hlam{a}$-module, 
$\overline{M} = \ftd M$ for some $M$ with 
$\epsi(M) =0$. Take $P = \Ind_{n,1}^{n+1} (q^i K_n \boxtimes q^i J_m)$.
Then $\f_{q^i J_m}(\overline{M}) = \f_P (M)
= \ft^{n+1}M = \ft \overline{M}$
is irreducible, as desired.
\end{proof}
This has the following simple consequence:

\begin{cor}
\label{cor-multfree}  \cite{GV}
\begin{enumerate}
\item[(i)]  Let $M$ be a simple $\Hlam{n}$-module, $\lambda$ cyclotomic.
Then the cosocle of $\Ind_{\Hlam{n}}^{\Hlam{n+1}} M$ and the
socle of $\Res_{\Hlam{n}}^{\Hlam{n+1}} M$ are multiplicity free.
\item[(ii)]  If $M$ is a simple $\Haffn$-module, 
$\soc \Res_{\Haff{n-1}}^{\Haff{n}} M$ is multiplicity free.
\end{enumerate}
\end{cor}
\begin{cor} If $M$ is irreducible, $\est_iM$ and $\fst_i M$ are indecomposable.
\end{cor}
\begin{prop}

\begin{enumerate}
\item[(i)]  $\soc (\esti M) \isom \cosoc(\esti M) $, for $M \in \RHlam$.
\item[(ii)] $\soc(\fsti M) \isom \cosoc (\fsti M)$, for $M \in \RHlam$,
$\lambda$ cyclotomic.
\end{enumerate}
\end{prop}

%
\begin{warning}
However, $\soc \Ind_{\Haff{n-1}\otimes \Haff{1}}^{\Haff{n}}
(M \boxtimes q^i J_1) \neq \cosoc \Ind_{\Haff{n-1}\otimes \Haff{1}}^{\Haff{n}}
(M \boxtimes q^i J_1)$.
For example, if $n=2$, and $M = q^{i-1} J_1$, they are clearly different.
\end{warning}

We finish the section with some alternate characterizations of $\epsi$
and $\nui$.

\begin{thm}
\label{thm-epsilon}
Let $N$ be an irreducible module in $\Rlam$.
Then
\begin{enumerate}
\item[(i)]  $\epsi (N)$ is the maximal size of a 
Jordan block for $R[[X_n - q^i]]$ on $N$.
\item[(ii)] In $K(\Rlam)$, we have 
$\esti N = \epsi (N) \cdot \et (N) +
\sum a_\alpha M_\alpha$,
where $M_\alpha$ are irreducible modules with
$\epsi (M_\alpha) < \epsi (N) -1$.
\item[(iii)]  $\epsi(N) = \dim \Hom(\esti N, \esti N)$.
\end{enumerate}
\end{thm}
\begin{proof}
We may suppose that $N = \ft^n M =
\cosoc \Ind (M \boxtimes q^i K_n)$,
where $M$ is an irreducible $\Hlam{a}$-module
and $\epsi (M) =0$.
By the Mackey formula and ``shuffle lemma,''
if $P \twoheadrightarrow q^i K_n$ is some module
covering $q^i K_n$ we have an exact sequence
$$
0 \to \Ind_{a, n-1,1}^{a+n-1,1}
\Res_{a, n-1,1}^{a, n}(M \boxtimes P)
\to \Res_{a, n-1,1}^{a+ n}
\Ind_{a, n}^{a+n}(M \boxtimes P) \to Q \to 0
$$
where $X_n$ has no $q^i$-eigenspace on $Q$, i.e.
$$
\Hom_{\Haff{a+n-1} \otimes \Haff{1}}
(\Haff{a+n-1} \boxtimes q^i J_1, Q) = 0.
$$
Hence taking the $q^i$-eigenspace of $X_n$, we get
$$
\esti \Ind_{a, n}^{a+n}(M \boxtimes P)
= \Ind_{a, n-1}^{a+n-1}(M \boxtimes \esti P)
$$
and that 
$\epsit (\Ind(M \boxtimes P) ) = \epsit(P)$,
where we denote the maximal size of a Jordan block of
$R[[X_n - q^i]]$ on a module $N$ by $\epsit (N)$.
We now apply this for $P = q^i K_n$.
Write, for $m \le n$
$$
0 \to S_m \to
\Ind(M \boxtimes q^i K_m) \to \ft^m M \to 0
$$
for the exact sequence defining $\ft^m M$.
We prove (ii).
Apply $\esti$ to this exact sequence, to get
$$0 \to \esti S_n \to
\Ind(M \boxtimes \esti(q^i K_n)) \to \esti \ft^n M \to 0.
$$
Now, in $K(\Rlam)$, $\esti(q^i K_n)$ is a multiple
of $q^i K_{n-1}$, as this is the unique irreducible
module with this central character.  Comparing
dimensions, we see that $\esti(q^i K_n) =
n \cdot q^i K_{n-1}$,
and so in $K(\Rep_{a+n}^\aff)$ we have
$$\esti \ft^n M =
n \cdot \Ind(M \boxtimes q^i K_{n-1}) - \esti S_n
= n \ft^{n-1} M + n S_{n-1} - \esti S_n.$$
By proposition \ref{prop-epsilon}, (ii), every subquotient
$A$ of $S_n$ has $\epsi (A) < m$, and
so the right hand side of the proposition is
of the desired form.
Finally, as $\ft^n M$ is in 
$K(\RHlam)$, it follows the non-zero
terms in this expression are also (after
cancelling).

We now prove (i). As $\Ind(M \boxtimes q^i K_n)
\twoheadrightarrow \ft^n M$, ({$\ast$}) implies
$$ \epsit (\ft^n M) \le \epsit(q^i K_n).$$
But by proposition \ref{prop-epsilon}, (ii) and its
proof, the $\Hlam{a} \otimes \Haff{n}$-module
$M \boxtimes q^i K_n$ does  occur as a
submodule of $\Res(\ft^n M)$ (in fact,
with multiplicity one).  It follows that
$$ \epsit (\ft^n M) \ge \epsit(q^i K_n).$$
It remains to show $\epsit(q^i K_n) = n =
\epsi(q^i K_n)$.

We show inductively that if $\alpha_m \neq 0$
in the exact sequence
$$0 \to \Hom(q^i J_{m-1}, q^i K_n)
\to \Hom(q^i J_{m}, q^i K_n)
\xrightarrow{\alpha_m} \Hom(q^i J_{1}, q^i K_n)$$
then $\Hom(q^i J_{m}, q^i K_n) = m \cdot q^i K_{n-1}$.
in the Grothendieck group of $\Haff{n-1}$ modules.
As $\epsit(q^i K_n)$ is the smallest
integer $m$ for which $\alpha_{m+1} = 0$,
and $\Hom(q^i J_{\epsit(q^i K_n)}, q^i K_n)
= \esti (q^i K_n)$, this will finish us up.
But $\Hom(q^i J_{1}, q^i K_n)$ is 
$\et(q^i K_n)$, which is $q^i K_{n-1}$ as it
is irreducible and non-zero.
Our induction starts with $m=0$, where it
is clear by these remarks. Likewise,
if $m >0$, then if $\alpha_m$ is non-zero it
must be a surjection onto $q^i K_{n-1}$,
hence the result.

We now show that (i) and (ii) imply (iii). 
By (i), the operators of multiplication by $(X_n - q^i)^k$
are non-zero for $0 \leq k < \epsi(N)$, and zero for $k = \epsi(N)$.
Hence they are linearly independent, and 
$\dim\Hom( \esti N, \esti N) \geq \epsi(N)$. For the reverse 
inequality, observe that for any module $X$ with cosocle $\et N$
$$ \dim\Hom(X,\esti N) \leq \epsi(N), $$
 as one sees by filtering $\esti N$ so it has semisimple quotients and applying (ii).
Now take $X = \esti N$.
\end{proof}

Let us write, for an irreducible $\Haffn$-module $N$ 
$$ 
        \epsihat(N) = \epsi(\sigma^* N) 
$$ 
where $\sigma\cc\Haffn \to \Haffn$ $T_i \mapsto -(T_{n-i} +1 -q)$, 
$ X_i \mapsto X_{n+1-i}$ is the 
{\it diagram automorphism\/} of $\Haffn$. 
Then $\epsihat(N)$ is the maximum, as $s \in (R^\times)^n$ varies 
with $N_s \neq 0$, of the length of the chain of $q^i$s beginning 
$s$.  (This is immediate as $\epsi(N)$ is the maximal 
length of the chain of $q^i$s ending $s$ with $N_s \neq 0$.) 
The previous theorem tells us $\epsihat(N)$ is also the maximal 
size of a Jordan block for $R[[X_1 - q^i]]$ on $N$; hence 
\begin{cor} 
\label{cor-epscheck}
If $N$ is an irreducible $\Haffn$-module, then 
$$      \prlam N \neq 0 \quad \iff \quad 
        \epsihat(N) \le \lambda_i, \qquad \text{for all $i \in \mu_q$.} 
$$ 
\end{cor} 
We now give descriptions of $\nui$ parallel to that of $\epsi$. 
The proof of the following theorem will not be completed until 
section \ref{sec-groth}.

\begin{thm}
\label{thm-nu}
Let $N$ be an irreducible module in $\Rlam$,
where $\lambda$ is cyclotomic. Then
\begin{enumerate}
\item[(i)]   $\nu_i(N)$ is the smallest integer $m$ for
which $\fsti N = \prlam \Ind (N \boxtimes q^i J_m)$.

\item[(ii)]   In the Grothendieck group $K(\Rlam)$,
$$ \fsti N = \nu_i(N) \cdot \ft (N) +
\sum a_\alpha^{\prime} M_\alpha^{\prime}$$
where $M_\alpha$ are irreducible modules with
$\epsi(M_\alpha') < \epsi(N) +1$.
\item[(iii)]  $\nui(N) = \dim \Hom(\fsti N, \fsti N)$.
\end{enumerate}
\end{thm}
\begin{proof}
We first show (i) and (ii) are equivalent.
Suppose the surjective map 
$$\prlam \Ind (N \boxtimes q^i J_m) \twoheadrightarrow
\prlam \Ind (N \boxtimes q^i J_{m-1})
$$
is {\it not\/} an isomorphism.  
We claim that it follows that 
$$\prlam \Ind (N \boxtimes q^i J_m) = m \cdot \ft(N) +
\text{ smaller terms},
$$
(where we mean equality in the Grothendieck group,
and {\it smaller terms\/} means a sum of modules
$A$ for which $\epsi(A) < \epsi ( \ft(N) )$.
The proof is by induction.
We have an exact sequence
$$\prlam \Ind (N \boxtimes q^i J_1) \xrightarrow{\alpha_m}
\prlam \Ind (N \boxtimes q^i J_m) \xrightarrow{\gamma_m} 
 \prlam \Ind (N \boxtimes q^i J_{m-1}) \to 0.
$$
By assumption, $\alpha_m \neq 0$, so the image
of $\alpha_m$ contains a copy of the irreducible module
$\cosoc \prlam \Ind (N \boxtimes q^i J_1) = \ft(N)$.
As $\prlam \Ind (N \boxtimes q^i J_1) = \ft(N)  +
\text{ smaller terms}$, the induction is complete if we show that
$\gamma_m$ is an isomorphism implies $\gamma_{m+k}$ is also, 
for $k \ge 0$. 
But this is clear---factor  
$\alpha_{m+k} = \theta_m^k \circ \alpha_m$, where 
$\theta_m^k$ is the map induced from $x^k\cc J_m \to J_{m+k}$. 
But now $\alpha_m = 0$ and so $\alpha_{m+k} =0$ and 
hence $\gamma_{m+k}$ is an isomorphism. 
 
Write $\widetilde{\varphi}_i(N)$ for the smallest integer $m$ 
for which $\fsti N = \pr \Ind(N \boxtimes q^i J_m)$. 
We have just shown that for $m > \widetilde{\varphi}_i(N)$,
$N \boxtimes q^i J_{\widetilde{\varphi}_i N}$ is a submodule of
$\Res \pr \Ind(N \boxtimes q^i J_m)$. It follows that for such $m$,
$ R^{\widetilde{\varphi}_i N} = \Hom( N \boxtimes q^iJ_m, N \boxtimes q^i J_{\widetilde{\varphi}_i N} )$
is a submodule of 
\begin{align*} 
\Hom (\Ind(N \boxtimes q^i J_m), 
        \pr \Ind(N \boxtimes q^i J_m) ) = 
\Hom(\fsti N, \fsti N),
\end{align*} 
and so $\dim \Hom(\fsti N, \fsti N) \geq \widetilde{\varphi}_i(N)$.
On the other hand, as in the proof of \ref{thm-epsilon}, for any module 
$X$ with cosocle $\ft N$, we have 
$$ \dim \Hom(X,\fsti N ) \leq \widetilde{\varphi}_i(N).$$
Taking $X = \fsti N$ we see 
(iii) is equivalent 
to (i) and (ii). 

Now let us show ${\widetilde{\varphi}_i (N)} \ge \nui(N)$.
But this is immediate from Frobenius reciprocity:  any 
embedding $R[X_1]/(X_1 - q^i)^m \inj q^i K_m$ 
induces a map $\Ind(N \boxtimes q^i J_{\nui(N)}) 
\to \Res \pr \Ind(N \boxtimes q^i K_{\nui(N)}),$ 
which when restricted to $\Haffn \otimes \Haff{1}$ induces an 
injection $N \boxtimes q^i J_{\nui(N)} \inj \Res \ft^{\nui(N)} N$. 
It follows that ${\widetilde{\varphi}_i (N)} \ge \nui(N)$. 
The reverse inequality will be proved in section \ref{sec-groth}. 
\end{proof}


\section{Serre relations}
\label{sec-serre}

In this section we study  the relation of 
$\epsi (\widetilde{\f}_{i q^{\pm 1}}N)$ to $\epsi (N)$.
To simplify notation we will write $i=1$, so
$i q = q$; 
this has no effect other than to make the notation readable.
We will also write $K_m$ rather than $q^0K_m$, and similarly for $J_m$,
and we shall write $q^i$ rather than $q^iJ_1$ when this causes no confusion.

The case $q=-1$ is different from the case $q \neq -1$;
for most of the section we will assume $q \neq -1$.
The main results of the section are proposition \ref{prop-handf},
which holds for $q \neq -1$, and its counterpart
proposition \ref{prop-q-1}, which covers the case $q=-1$.

\noindent \underline{$q \neq -1$}   
First, suppose $q \neq -1$.
Then we have a short exact sequence of
$\Haff{2}$-modules 
$$0 \to \triv \to \Ind (q \boxtimes 1) \to \St \to 0$$
which does not split;
here $\triv$ and $\St$ are the one
dimensional modules with spectra $(1q)$ and $(q1)$ respectively.

Hence for $m \ge 1$ we have a short exact sequence ({$**$})
$$0 \to \Ind (N \boxtimes \triv \boxtimes K_{m-1}) \to
 \Ind(N \boxtimes q \boxtimes K_m) \to
 \Ind(N \boxtimes \St\boxtimes K_{m-1}) \to 0.$$

\begin{prop}
\label{prop-f1q}
Suppose $q \neq -1$.
Then for all irreducible $N \in \Rep \Hlam{a}$, and all $m \ge 1$
\begin{enumerate}
\item[(i)] $\Ind (\triv \boxtimes K_m) \isom 
\Ind (K_m \boxtimes \triv)$ is irreducible.
\item[(ii)] $\Ind(N \boxtimes  K_m\boxtimes \triv) $ is
indecomposable, with simple cosocle.
\item[(iii)]$ \Ind(\St \boxtimes  K_m) \isom
\Ind (K_m \boxtimes \St)$ is irreducible.
\item[(iv)]$\Ind(N \boxtimes  K_m\boxtimes \St)$
is indecomposable, with simple cosocle.
\end{enumerate} 
\end{prop}
\begin{proof}
We prove (i) and (ii); as (iii) and (iv) are similar
(or indeed are formal consequences, upon applying the obvious automorphism
of $\Haffn$).
By the classification of $\Haff{3}$-modules, when 
$q \neq -1$,
$$\Ind (\triv \boxtimes 1) \isom \Ind (1 \boxtimes \triv)$$
is irreducible; hence by transitivity of induction
$$\Ind (\triv \boxtimes 1 \boxtimes \cdots  \boxtimes 1) =
\Ind (1 \boxtimes \triv \boxtimes \cdots  \boxtimes 1) =
\cdots = \Ind (1 \boxtimes \cdots  \boxtimes 1 \boxtimes \triv).$$
But the Mackey formula implies that $\triv \boxtimes K_m$
occurs with multiplicity one in $\Res \Ind(\triv \boxtimes K_m)$,
hence the {\it useful observation\/} (lemma \ref{useful}) shows
$\Ind(K_m  \boxtimes \triv)$ is irreducible.

To prove (ii), observe that if $X \surj Y$ is a surjection,
then $\cosoc(X) \surj \cosoc(Y)$, so as $N$ is
a quotient of $\Ind(\widetilde{N} \boxtimes K_m)$ for
some $\widetilde{N}$ with $\varepsilon_1(N) =0$, it is
enough to show (ii) when $\varepsilon_1(N) =0$.
We play with the Mackey formula again as in theorem \ref{thm-epsilon},
so we will be terse.
Write $\Res = \Res_{a, m+2}^{a+m+2}$. Then the sequence
$$0 \to N \boxtimes \Ind (K_m  \boxtimes \triv) 
\to \Res \Ind (N \boxtimes K_m  \boxtimes \triv) \to U \to 0$$
splits, as the shuffle lemma shows the central character of 
$U$ differs from that of $N \boxtimes \Ind (K_m  \boxtimes \triv)$.
It follows, as in \ref{prop-epsilon}, that
$$\Res { } \cosoc \Ind (N \boxtimes K_m  \boxtimes \triv)=
N \boxtimes \Ind (K_m  \boxtimes \triv) \oplus \overline{U},$$
%
where $\overline{U}$ is some quotient of $U$, and that
\begin{align*}
\Hom(\cosoc \Ind (N \boxtimes K_m  \boxtimes \triv),
\cosoc \Ind (N \boxtimes K_m  \boxtimes \triv)) &= \\
\Hom(N \boxtimes \Ind (K_m  \boxtimes \triv), N \boxtimes \Ind (K_m  
\boxtimes \triv)) &= R,
\end{align*}
the last equality as $\Ind(K_m  \boxtimes \triv)$ is irreducible by 
(i).
\end{proof}

For $N \in \B$, define $\fq{i }{ q\!i} (N) = \prlam  \cosoc 
\Ind(N \boxtimes {\bf 1}_{(i q\!i)})$, and 
$\fq{q\!i }{ i}(N) =  \prlam \cosoc$  $ \Ind(N \boxtimes \st{q\!i}{ i})$.
As a corollary of the proposition, $\fq{q\!i }{ i}N $ and $\fq{i}{ q\!i}N$
are irreducible modules, if $q \neq -1$.
(If $q=-1$ this is no longer true.)

\begin{lemma}
\begin{enumerate}
\item[(i)] If $A$ is an irreducible constituent of 
$\Ind(N \boxtimes \triv \boxtimes K_{m-1})$, then
$\varepsilon_1(A) \le \varepsilon_1(N) + m -1$.
\item[(ii)] If $A$ is an irreducible constituent of 
$\Ind(N \boxtimes K_{m-1} \boxtimes \St)$, then
$\varepsilon_1(A) \le \varepsilon_1(N) + m$. Moreover, if 
$\varepsilon_1(A) = \varepsilon_1(N) + m$,
then $A$ is an irreducible constituent  of 
$\Ind(\widetilde{\f}_1^{m-1} N \boxtimes \St)$.
	\item[(iii)] If $Q$ is an irreducible quotient of 
$\Ind(N \boxtimes K_{m-1} \boxtimes \St)$, then $\varepsilon_1(Q)
= \varepsilon_1(N) + m$.
\end{enumerate}
\end{lemma}
\begin{proof}
(i) and the first part of (ii) are immediate from the shuffle lemma.
As $\Ind(N \boxtimes K_{m-1}) = \widetilde{\f}_1^{m-1} N
+ \sum a_\alpha M_\alpha$, 
where $\varepsilon_1(M_\alpha) < \varepsilon_1(N) + m -1$
by proposition \ref{prop-epsilon}, the last part of (ii) follows by another
application of the shuffle lemma.

Finally, to prove (iii), it is enough to show $\varepsilon_1(Q)
\ge  \varepsilon_1(N) + m$.
Further, we may again assume  $\varepsilon_1(N) =0$, as the cosocle of
$\Ind(\widetilde{\e}_1^{\varepsilon_1(N)}N \boxtimes K_{
\varepsilon_1(N) +m-1}  \boxtimes \St)$ surjects onto the
cosocle of $\Ind(N \boxtimes K_{m-1} \boxtimes \St)$.
Now, if $Q$ is an irreducible quotient as in (iii),
Frobenius reciprocity gives a non-zero homomorphism
$N \boxtimes \Ind(K_{m-1} \boxtimes \St) \to Q$;
hence $\varepsilon_1(Q) \ge \varepsilon_1(N \boxtimes \Ind(K_{m-1}
 \boxtimes \St))$.
But $\Ind(K_{m-1}
 \boxtimes \St) = \Ind(\St \boxtimes K_{m-1})$ by the previous
proposition, and so $\varepsilon_1(Q) \ge m$.
\end{proof}

We have shown that
$$\Ind(N \boxtimes \triv \boxtimes K_{m-1}) = 
\fq{1}{q} \widetilde{\f}_1^{m-1} N + \text{ smaller terms,}$$
and
$$\Ind(N \boxtimes \St \boxtimes K_{m-1}) = \cosoc + \text{ smaller terms,}$$
where all the terms in the cosocle have $\varepsilon_1 =
\varepsilon_1(N)  +m$. As
$\Ind(N \boxtimes q\boxtimes K_m)$ surjects onto
$\widetilde{\f}_1^m \widetilde{\f}_q N$ for $m \ge 0$, and
we have a filtration of $\Ind(N \boxtimes q\boxtimes K_m)$
as in (${**}$), we have proved most of

\begin{prop}
Precisely one of the following alternatives hold
\begin{enumerate}
\item[(i)] For all $m \ge 1$, $\widetilde{\f}_1^m \widetilde{\f}_q N
= \fq{1}{q} \widetilde{\f}_1^{m-1} N$, and for all
$m \ge 0$ 
$$\varepsilon_1(\widetilde{\f}_1^m \widetilde{\f}_q N) =
m-1+ \varepsilon_1(N), \qquad \text{ or}$$
\item[(ii)] For all $m \ge 1$, $\widetilde{\f}_1^m \widetilde{\f}_q N$
is a summand of $\cosoc \Ind(\widetilde{\f}_1^{m-1} N \boxtimes \St)$,
and for all $m \ge 0$,
$$\varepsilon_1(\widetilde{\f}_1^m \widetilde{\f}_q N) =
m +\varepsilon_1(N).$$
\end{enumerate}
\end{prop}
\begin{proof}
We have proved everything for $m \ge 1$; to finish we must only
observe that $\varepsilon_1(\widetilde{\f}_1 \widetilde{\f}_q N)
= \varepsilon_1(\widetilde{\f}_q N) +1$, and so we have
the assertions for $m=0$ also.
\end{proof}

Now let $N$ be an irreducible $\Hlam{n}$-module, with $\lambda$
cyclotomic. 
Let us agree to write $\ft$ for the {\it affine\/} crystal operator
$\cosoc \Ind( \bullet \boxtimes q^i)$, and write 
$\pr \ft$ for the {\it cyclotomic\/} crystal operator.

\begin{prop}
\label{prop-handf}
Let $N$ be an irreducible $\Hlam{n}$-module, and suppose
$\prlam \widetilde{\f}_q N \neq 0$.
Then precisely one of the following holds:
\begin{enumerate}
\item[(i)]$ \varepsilon_1(\widetilde{\f}_q N) =
\varepsilon_1(N) -1$, and $\varphi_1(\widetilde{\f}_q N) =
\varphi_1(N), \qquad$ or
\item[(ii)] $\varepsilon_1(\widetilde{\f}_q N) =
\varepsilon_1(N),$ and  $\varphi_1(\widetilde{\f}_q N) =
\varphi_1(N) +1$.
\end{enumerate}
\end{prop}
\begin{proof}
Write $\varphi = \varphi_1(N).$ First suppose 
$\varepsilon_1(\widetilde{\f}_q N) = \varepsilon_1(N) -1$.
Then the surjection $\Ind(1 q) \to \triv$ induces a surjection,
hence isomorphism
\begin{align*}
\widetilde{\f}_q \widetilde{\f}_1^m N =
\cosoc \Ind(N \boxtimes K_m \boxtimes q) \surj
&\cosoc \Ind(N \boxtimes K_{m-1} \boxtimes \triv) \\ 
&= \fq{1}{q} \widetilde{\f}_1^{m-1} N = \widetilde{\f}_1^m
\widetilde{\f}_q  N
\end{align*}
for all $m \ge 1$.
It follows that $\prlam \widetilde{\f}_1^{\varphi+1} \widetilde{\f}_q  N
= \prlam \widetilde{\f}_q  \widetilde{\f}_1^{\varphi+1} N = 0$.
To show $\prlam \widetilde{\f}_1^{\varphi} \widetilde{\f}_q  N \neq 0$,
we compute $\ech{\alpha}(\widetilde{\f}_1^m 
\widetilde{\f}_q  N)$ for all $\alpha \in \mu_q$ and
invoke corollary \ref{cor-epscheck}.

The shuffle lemma shows 
$\echone(\Ind(\widetilde{\f}_1^m  N\boxtimes q)) \le
\echone(\widetilde{\f}_1^m N)$,
hence the definition of $\varphi$
and corollary \ref{cor-epscheck} shows $\echone(\widetilde{\f}_1^\varphi
\widetilde{\f}_q  N) \le \lambda_1$.

Similarly, the shuffle lemma shows that for all $m$,
$\ech{q} \Ind(\widetilde{\f}_q  N \boxtimes  K_m)
\le \ech{q}(\widetilde{\f}_q N)$ 
and $\ech{q}(\widetilde{\f}_q N) \le \lambda_q$
by corollary \ref{cor-epscheck} and the assumption 
$\prlam \widetilde{\f}_q N \neq 0$.
We clearly have $\ech{\alpha}(\widetilde{\f}_1^m
\widetilde{\f}_q  N) \le \ech{\alpha}(N) \le \lambda_\alpha$
for $\alpha \notin \{1,q\}$, and all $m \ge 0$.
Thus another application of corollary \ref{cor-epscheck} shows 
$\prlam(\widetilde{\f}_1^\varphi \widetilde{\f}_q  N) \neq 0$,
so $\varphi_1(\widetilde{\f}_q  N) = \varphi$.

Now suppose $\varepsilon_1(\widetilde{\f}_q N) = \varepsilon_1(N)$.
The argument of the last paragraph applies equally well here,
so to show $\prlam(\widetilde{\f}_1^{\varphi +1} \widetilde{\f}_q  N) \neq 0$
we must only show $\echone(\widetilde{\f}_1^m \widetilde{\f}_q  N)
\le \echone(\widetilde{\f}_1^{m-1} N) $
for all $m \ge 1$. But this is immediate from the shuffle lemma, as
in this case $\widetilde{\f}_1^m \widetilde{\f}_q  N$ is a quotient
of $\Ind(\widetilde{\f}_1^{m-1} N \boxtimes \St)$.

Finally, observe that Frobenius reciprocity implies that 
$\widetilde{\f}_1^{m-1} N \boxtimes \St$ is contained in
$\Res(\widetilde{\f}_1^m \widetilde{\f}_q  N)$,
so if $\echone(\widetilde{\f}_1^{m-1} N) > \lambda_1$,
then $\echone(\widetilde{\f}_1^m \widetilde{\f}_q  N) > \lambda_1$ also.
It follows that $\prlam \widetilde{\f}_1^{\varphi + 2} \widetilde{\f}_q  N 
=0$, and so $\varphi_1(\widetilde{\f}_q  N) = \varphi +1$ in this case.
\end{proof}

\noindent \underline{$q = -1$}   We now suppose $q = -1$.
The analogue of proposition \ref{prop-handf} is

\begin{prop}
\label{prop-q-1}
Let $N$ be an irreducible $\Hlam{n}$-module, and suppose
$\prlam \widetilde{\f}_q N \neq 0$.
Then precisely one of the following holds:
\begin{enumerate}
\item[(i)]$ \varepsilon_1(\widetilde{\f}_q N) =
\varepsilon_1(N) -2$, and $\varphi_1(\widetilde{\f}_q N) =
\varphi_1(N), \qquad$ or
\item[(ii)]$ \varepsilon_1(\widetilde{\f}_q N) =
\varepsilon_1(N) -1$, and $\varphi_1(\widetilde{\f}_q N) =
\varphi_1(N) +1, \qquad$ or
\item[(iii)] $\varepsilon_1(\widetilde{\f}_q N) =
\varepsilon_1(N),$ and  $\varphi_1(\widetilde{\f}_q N) =
\varphi_1(N) +2$.
\end{enumerate}
\end{prop}
\begin{proof}
We merely sketch the necessary modifications in the definitions
needed to prove this.
Recall that there are 3 irreducible representations of
$\Haff{3}$ with central character $S_3\cdot (1 1 q)$.
Denote them $\alpha, \overline{\alpha}, \gamma$
with $\ch \alpha = 2(1 1 q)$, $\ch \overline{\alpha} = 2(q 1 1)$,
$\ch \gamma = (1 q 1)$.
The representations $\Ind(1 \boxtimes \alpha)$ and 
$\Ind(1 \boxtimes \gamma)$ are irreducible, so if we define
\begin{gather*}
\widetilde{\f}_{11q} (M) = \prlam \cosoc \Ind(M \boxtimes \alpha) \\
\widetilde{\f}_{1q1} (M) = \prlam \cosoc \Ind(M \boxtimes \gamma ) \\
\widetilde{\f}_{q11} (M) = \prlam \cosoc \Ind(M \boxtimes
\overline{\alpha})
\end{gather*}
then the analogue of proposition \ref{prop-f1q} is that
each of these operators takes irreducible modules to irreducible
modules (or zero).
Using the exact sequence $0 \to \gamma \to
\Ind(1 \boxtimes \triv) \to \alpha \to 0$, we see that for
$m \ge 2$ there is a 4 step filtration of $\Ind(N \boxtimes q \boxtimes K_m)$
$$0 = F_0 \subset F_1 \subset F_2 \subset F_3 \subset F_4 =
\Ind(N \boxtimes q \boxtimes K_m)$$ 
where
\begin{gather*}
F_1 = \Ind(N \boxtimes \Ind(\alpha \boxtimes K_{m-2})),
\qquad F_2/F_1 = F_3/F_2 = \Ind(N \boxtimes \Ind(\gamma \boxtimes K_{m-2})),\\
F_4/F_3 = \Ind(N \boxtimes \Ind(\overline{\alpha} \boxtimes K_{m-2})).
\end{gather*}
Arguing as before, one sees that the cosocle of each
of these subquotients consists of a sum of terms with fixed
$\varepsilon_1$: 
respectively $\varepsilon_1 -m$ is $\varepsilon_1(N) -2$,
 $\varepsilon_1(N) -1$, and  $\varepsilon_1(N)$.
Hence $\widetilde{\f}_1^m \widetilde{\f}_q N$, which is a quotient
of $\Ind(N \boxtimes q \boxtimes K_m)$, is a quotient
of precisely one of those three subquotients (for $m \ge 2$).
The rest of the proof is as before.
\end{proof}

\section{Shapovalov Form}
\label{sec-shap}
Let $A$ be a finite dimensional algebra over $R$, and 
$\Proj(A)$ be the category of finite dimensional projective
$A$-modules.
As $\Hom_A(P, - )$ is an exact functor on $\Rep A$ if $P$ is
projective, there is a well defined bilinear form
$$ ( \quad, \quad)\cc  K(\Proj A) \otimes_\Z K(\Rep A) \to \Z
$$
given by $([P], [M]) = \dim \Hom_A(P,M)$.

If $M$ is a simple $A$-module, we write $\P{M}$ for the
projective cover of $M$, so $\P{M} $ is an indecomposable
projective with $\cosoc(\P{M}) = M$.
Then $(\P{M}, N)$ is the multiplicity of $M$ in the
composition series of $N$, i.e.~%
$(\P{M}, - ) = \delta_M\cc  K(A\text{-mod}) \to \Z$,
and $( \quad, \quad)$ defines an isomorphism 
$K(\Proj A) \xrightarrow{~} K(\Rep A)^*$.

Now consider the natural map
$$K(\Proj A) \to K(\Rep A).
$$
This is an injection if and only if it becomes an isomorphism
after tensoring $\otimes_\Z \Q$, and this happens precisely when
the bilinear form in non-degenerate when restricted to
$$K(\Proj A) \otimes_\Z K(\Proj A) \to \Z.
$$
We now show this happens for the cyclotomic Hecke algebra.
In this section $\Hlam{n}$ denotes the cyclotomic
Hecke algebra only, i.e.~%
$\lambda\cc  \mu_q \to \Z_+ $ is a function with 
$\sum \lambda_i < \infty$.
\begin{thm}
\label{thm-shap}
The above pairing $(P, P') \mapsto \dim \Hom(P,P')$
$$ ( \quad, \quad)\cc  K(\RHlam)^* \otimes_\Z K(\RHlam)^* \to \Z
$$
is a non-degenerate symmetric form.
\end{thm}
We call this form, or the induced form
$$ K(\RHlam) \otimes_\Z K(\RHlam)\to \Q
$$
the Shapovalov form.

Observe that as $\esti$ and $\fsti$ carry projective modules
to projective modules, they act on $K(\Proj \Hlam{a})$.
Clearly the action is just the transpose of the action of 
$\esti$ and $\fsti$ on $K(\Rep \Hlam{a})$.
Hence the following lemma is an immediate consequence of
the results of section \ref{sec-crystal}. (We reassure the reader
that the results of this section are not used in the proof of
theorem \ref{thm-nu}).
\begin{lemma}
\label{lemma-shap1}
Let $N$ be an irreducible $\Hlam{a}$-module, and set
$\varepsilon = \epsi(N)$, $\nu = \nui(N)$. Then
\begin{enumerate}
\item[(i)] $\esti^{(\varepsilon)} \P{N} =
\binom{\varepsilon + \nu}{\varepsilon} \P{\et^\varepsilon N}
+ \sum_{Q: \epsi(Q) > \varepsilon} a_Q \P{\et^\varepsilon Q}$,
for some $a_Q \in \N$.
\item[(ii)] If $n > \varepsilon$, $\esti^{(n)} \P{N} =
\sum_{Q: \epsi(Q) \ge n > \varepsilon} \alpha_Q \P{\et^n Q}$ 
for some $\alpha_Q \in \N$.
\end{enumerate}
\end{lemma}
\begin{proof}
If $S$ is an irreducible $\Hlam{a-\varepsilon}$-module, then
$\dim \Hom(\esti^{(\varepsilon)} \P{N}, S) =$
$\dim$ \linebreak[4] $\Hom(\P{N}, \fsti^{(\varepsilon)} S)$ is
clearly zero unless $\nui(S) \ge \varepsilon$.
So we may assume $S = \et^\varepsilon Q$, where 
$Q$ is an irreducible $\Hlam{a}$-module, and then
$$a_Q = \dim \Hom(\esti^{(\varepsilon)} \P{N}, \et^\varepsilon Q)
= \dim \Hom(\P{N}, \fsti^{(\varepsilon)} \et^\varepsilon Q).
$$
But $\fsti^{(\varepsilon)} \et^\varepsilon Q 
= \binom{\nui(Q) + \varepsilon}{\varepsilon }Q +
\sum_{\epsi(M) < \epsi(Q)} \alpha_M M$, by proposition \ref{prop-epsilon}.

Now suppose that $\epsi(Q) = \varepsilon = \epsi(N)$.
Then the terms $M$ in the above sum have $\epsi(M) < \varepsilon$,
so none are isomorphic to $N$.
It follows that if $\epsi(Q) = \varepsilon$, $a_Q = 0$ unless
$Q= N$, and $a_N  = \binom{\nu + \varepsilon}{\varepsilon}$.
The lemma is immediate.
\end{proof}
We now prove the theorem.
\begin{proof}
We may suppose that projective covers of irreducible
modules in $K(\Rep \Hlam{a'})$ are linearly independent
for $a' < a$, and that $a > 0$.
Suppose we have a relation
$$\sum c_M \P{M} = 0
$$
in   $K(\Rep \Hlam{a})$ with not all the $c_M$ equal
to zero.
Choose an $i \in \mu_q$ and a simple module $N$
such that $c_N \neq 0$, and $\varepsilon = \epsi(N)$ is
maximal among terms in this sum.
We may choose $i \in \mu_q$ so that $\varepsilon > 0$, 
as for a module $N$, $\epsi(N) = 0$ for all $i \in \mu_q$
implies that $\Res_{\Hlam{a-1}}^{\Hlam{a}} N = 0$,
which is absurd for $a >0$.

Now apply $\esti^{(\varepsilon)}$ to this sum.
By the lemma,
we get  an equality
$$\sum_{N: \epsi(N) = \varepsilon} \binom{\varepsilon + \nui(N)}{\varepsilon}
c_N \P{\et^\varepsilon N} + X =0
$$
in $K(\Rep \Hlam{a-\varepsilon})$, where $X$ is  a sum
of terms of the form $\P{\et^\varepsilon Q}$, with
$Q \in \Rep_q \Hlam{a}$ and $\epsi(Q) > \varepsilon$.
In particular, all the terms in the sum are
distinct projective modules in $\Rep \Hlam{a-\varepsilon}$.
By our inductive assumption these terms are linearly
independent, hence $X =  c_N = 0$.
This contradicts our choice of $c_N$, and shows
the Shapovalov form is non-degenerate.

It remains to show it is symmetric.
Again, induct on $a$.
Clearly the form is symmetric on $K(\Rep \Hlam{0})$,
and $(\fsti x, y) = (x, \esti y) = (\esti y, x) 
= (y, \fsti x)$ where we have used adjunction twice,
and the inductive hypothesis.
So the form is symmetric on the image of $\1$ under
the operators $\esti^{(n)}$ and $\fsti^{(n)}$.
We must merely show that this is everything.
This follows from the following two lemmas.
\end{proof}

\begin{lemma}
\label{lem-shap2}
Let $N$ be an irreducible $\Hlam{a}$-module, and set
$\varepsilon = \epsi(N)$, $\nu = \nui(N)$. Then
\begin{enumerate}
\item[(i)] $\fsti^{(\nu)} \P{N} =
 \binom{\varepsilon + \nu}{\varepsilon} \P{\ft^\nu N}
+ \sum_{Q: \nui(Q) > \nu} a_Q \P{\ft^\nu Q}$
for some coefficients $a_Q \in \N$.
\item[(ii)] If $n > \nu$, then 
 $\fsti^{(n)} \P{N} =
\sum_{Q: \nui(Q) \ge n > \nu} \alpha_Q \P{\ft^n Q}$
for some coefficients $\alpha_Q \in \N$.
\end{enumerate}
\end{lemma}
We omit the proof.
\begin{lemma}
\label{lem-shap3}
Every $\P{N} \in K(\Proj \Hlam{a})$ can be written as
a sum of monomial words in $\fsti^{(n)}$
with integer coefficients.
\end{lemma}
\begin{proof}
Again, we assume the result for $a' < a$, and suppose $a>0$.
Fix $i \in \mu_q$, and $r > 0$, and suppose
the result is true for all irreducible $N$ with 
$\nui(N) > r>0$.
As there are finitely many irreducible 
modules in $\RHlam$, our induction on $r$ starts successfully,
somewhere.
Let $M$ be irreducible, and $\epsi(M) = r$.
Apply the above lemma
to $N = \et^r M$, to get
$$\P{M} = \fsti^{(r)} \P{\et^r M} \; - \sum_{Q: \nui(Q) >r}
a_Q \P{\ft^r Q}
$$
which by induction is of the desired form.
Arguing in this way for each $i \in \mu_q$, 
we see the result is true except perhaps for modules
$N$ for which $\nui(N)=0$ for all $i \in \mu_q$.
Such a module would have $\Ind_{\Hlam{a}}^{\Hlam{a+1}} N =0$,
which is absurd.
Hence the result is true for all modules.
\end{proof}

\section{Relations in the Grothendieck group}
\label{sec-groth}
We show the operators $\esti, \fst_i, \hst_i$ acting on the Grothendieck
group of modules for the cyclotomic
and affine Hecke algebras satisfy the defining relations of $\sl$.
This is {\it not\/} true before passage to the Grothendieck group.
%
We know of several proofs of this result. 
Aside from the one given here, another would be to simply 
compute explicitly the action of $\esti$ and $\fsti$ on 
a suitable basis of $K(\Rlam)$. 
Such a basis is given by the  Specht modules, first defined 
by Ariki and Koike. 
 
In order to do this, one must (i) define these modules,  
(ii) show they span the Grothendieck group $K(\Rlam)$, and 
(iii) explicitly compute  $\esti$ and $\fsti$ on these modules. 
This requires quite some work (of a combinatorial nature) which 
can, however, be found in the literature. 
 
In the approach below, we derive the defining relations for 
$\sl$ from our general theory (which is built from an explicit 
study of the representations of $\Haffn, n \le 4$). 
As a consequence we rederive  and explain
these properties of Specht modules, obtaining a conceptual explanation 
for their combinatorics: it is just a realization of the crystal graph.

\begin{prop}
\label{prop-serre}
The operators $\esti\cc  K(\Rlam) \to K(\Rlam), \;$
$\esti\cc  K(\Raff) \to$ \linebreak[4]
 $K(\Raff)$ satisfy the Serre relations; i.e.~%
if $i, j \in \mu_q, \; i j^{-1} \neq q^{\pm 1},$
then as operators on the Grothendieck group
$$
\esti \e_j^{\ast} = \e_j^{\ast} \esti :   K(\Rlam) \to K(\Rlam).
$$
If $i j^{-1} = q^{\pm 1}$, and $q \neq q^{-1}$, then
$${\esti}^2  \e_j^{\ast}  + \e_j^{\ast} {\esti}^2 =  2 \esti \e_j^{\ast} \esti
 :   K(\Rlam) \to K(\Rlam).$$
If $i j^{-1} = q$ and  $q = q^{-1} = -1$ then 
$$ {e_i^*}^3e_j^* + 3{e_i^*}{e_j^*}{e_i^*}^2 =
   3{e_i^*}^2{e_j^*}{e_i^*} + {e_j^*}{e_i^*}^3 :  K(\Rlam) \to K(\Rlam).$$
\end{prop}
\begin{proof}
This reduces to checking the result on irreducible modules
for $\Haff{2}$ (if $i j^{-1} \neq q^{\pm 1}$),
$\Haff{3}$ (if $i j^{-1} = q^{\pm 1}$, and $q \neq q^{-1}$),
and $\Haff{4}$ (in the remaining case).
To see this, observe that as $\esti$ commutes with 
$\evst \cc  K(\Rlam) \to K(\Raff)$, and $\evst$ is injective, it is
enough to check that the $\esti$ satisfy the Serre relations
on $K(\Raff)$.
But, as observed in \ref{sec-divided}, $\esti$ is just the component 
$Id \otimes \delta_{q^i J_1}$ of $\Delta_{n,1}$ and $\Delta$ is
coassociative.  So it is enough to check the relations involving
a word of length $k$ in the $\esti$'s on $\Haff{k}$.

We must now check the Serre relations for the
irreducible $\Haff{n}$-modules, $n \le 4$.
These modules were listed in section \ref{sec-examples}; and the
relations follow by examining each module individually.
\end{proof}

\begin{rem}
The argument reducing the proposition to a case by case check 
sounds simpler if phrased in terms of 
$K(\Raff)^{\ast}$---algebras are easier to think with than
coalgebras.
It becomes even clearer when stated explicitly:
Let $M \in \Rep \Haff{n+2}.$  Consider $\esti \e_j^{\ast} (M)$.
This is $\varinjlim_{m, m^\prime} \Hom (
\Haff{n} \boxtimes q^i J_m \boxtimes q^j J_{m^\prime},
\Res_{\Haff{n} \otimes \Haff{1} \otimes \Haff{1}}^{\Haff{n+2}} M )$.
Factor $\Res_{n,1,1}^{n+2} M = \Res_{n,1,1}^{n,2} \Res_{n,2}^{n+2} M$.
$\Res_{n,2}^{n+2} M$ has a filtration with graded pieces
simple modules $N_\alpha \boxtimes \Gamma_\alpha$, with 
$N_\alpha$ a simple $\Haff{n}$-module, and $\Gamma_\alpha$
a simple $\Haff{2}$-module.
The image on $\esti \e_j^{\ast} (M)$ in the Grothendieck group
depends only on the value of the exact functor
\begin{multline*}
 \mathop{\varinjlim}\limits_{m, m^\prime}  \:
\Hom (\Haff{n} \boxtimes q^i J_m \boxtimes q^j J_{m^\prime},
\Res_{n,1,1}^{n+2}(\quad ) ) = \\
 Id \otimes \mathop{\varinjlim}\limits_{m, m^\prime}
\: \Hom_{R[X_{n+1}^{\pm 1}, X_{n+2}^{\pm 1}]} 
(q^i J_m \boxtimes q^j J_{m^\prime},
\Res_{1,1}^2 ( \quad ) )
\end{multline*}
on the  pieces $N_\alpha \boxtimes \Gamma_\alpha$,
and this functor depends only on $\Gamma_\alpha$.
As $\e_j^{\ast} \esti$ has a similar description, acting on the
pieces  $N_\alpha \boxtimes \Gamma_\alpha$ by
some other functor depending only on $\Gamma_\alpha$,
it suffices to check equality on $\Haff{2}$-modules.
And again, as these functors on $\Haff{2}$ are exact, it suffices
to check equality on simple modules.
\end{rem}
%
\begin{rem} 
In fact, the Serre relations are a formal consequence of the following 
(see section 3.3.3 of \cite{Kac}): 
(i) The relations $[\esti, \f_j^*] = \delta_{ij} \hsti$ 
and $[\hsti, \e_j^*] = c_{ij}  \e_j^*$, proved below,
(ii) the Shapovalov form of section \ref{sec-shap} is non-degenerate,
and (iii) the union $\bigcup \Rlam = \Raff$.
(This last condition ensures that the Serre relations hold not just as
operators on the irreducible module $K(\Rlam)$, but also
on the  module $K(\Raff)$.)
This is no simplification, as the proof we have chosen to give of
$[\hsti, \e_j^*] = c_{ij} \e_j^*$ essentially consists of verifying a
more precise form of the Serre relations.
\end{rem} 

\begin{cor}
The operators $\fsti \cc  K(\Rlam) \to K(\Rlam)$ satisfy
the Serre relations.
\end{cor}
\begin{proof}
It is enough to show they satisfy the Serre relations as 
maps from $K(\Rlam)_\Q \linebreak[3] \to K(\Rlam)_{\Q}$.
But the Shapovalov form of section \ref{sec-shap} is non-degenerate, and $\fst_i$
is adjoint to $\est_i$ with respect to this form.
As the $\est_i$ satisfy the Serre relations, so do the $\fst_i$.
\end{proof}

\begin{prop}
\label{prop-dividedpowers}
$n ! \edst = \esti^n, \qquad n ! \fdst = \fsti^n \cc  K(\Rlam) \to K(\Rlam)$
\end{prop}
\begin{proof}
As before, to check $n ! \edst = \esti^n$, it suffices to
check on simple $\Haff{n}$-modules.
Both left and right hand sides are zero on all simple
$\Haff{n}$-modules except $q^i K_n$, where both are $n !$.
This implies $n ! \fdst = \fsti^n$, by an argument similar
to the above one.
\end{proof}

We must now determine the relations between $\esti$ and $\f_j^*$. 
The Mackey formula for the cyclotomic Hecke algebra implies 
\begin{gather} 
\label{cyc-ind-re} 
 [\Res_{\Hlam{\bullet}}^{\Hlam{\bullet +1}}, 
\Ind_{\Hlam{\bullet}}^{\Hlam{\bullet +1}}](M) = 
(\sum \lambda_i) M. 
\end{gather} 
As $\esti$ is a refined version of restriction, and  
$\fsti$ is a refined variant of induction, 
the next theorem should be regarded as a sharpening 
of \eqref{cyc-ind-re}. 
The proof will occupy the rest of this section.

Define for $N$ an irreducible $\Hlam{n}$-module 
$$\hsti(N) = (\nui(N) -\epsi(N)) N$$ 
and more generally  define $\binom{\hsti}{k} (N) = 
\binom{\nui(N) -\epsi(N)}{k} \cdot N$, 
so that $\binom{\hsti}{k}\cc K(\Rlam) \to K(\Rlam)$. 
\begin{thm} 
\label{thm-[e,f]}
$[\esti, \f_j^*] = \delta_{ij} \: \hsti \cc K(\Rlam) \to K(\Rlam)$. 
\end{thm} 
We begin by showing 
\begin{prop} 
\label{prop-hM}
Suppose $M$ is an irreducible $\Hlam{n}$-module. 
Then $[\esti, \fsti](M)$ is a multiple of $M$, and 
$[\esti, \f_j^*](M) =0$ if $i \neq j$. 
\end{prop} 
\begin{proof} 
For $m \gg 0$ we have a surjection 
$$\Ind(M \boxtimes q^j J_m) \surj \f_j^* M \to 0.$$ 
Apply $\pr \esti$. As 
$\esti$ is exact, and $\pr$ is right exact, we still get a  
surjection 
$$ 
        \pr \esti \Ind(M \boxtimes q^j J_m) \surj \esti \f_j^* M \to 0. 
$$ 
But by the Mackey formula, we have an exact sequence 
$$ 
        0 \to \delta_{ij} m M \to 
        \esti \Ind(M \boxtimes q^j J_m) \to 
        \Ind(\esti M \boxtimes q^j J_m) \to 0, 
$$ 
and hence, as $M$ is irreducible, an exact sequence 
$$ 
        0 \to \delta_{ij} m' M \to 
        \pr \esti \Ind(M \boxtimes q^j J_m) \to \f_j^* \esti M \to 0 
$$ 
for some $m' \le m$, if $m \gg 0$. Hence 
\begin{gather} 
\label{gg} 
        \delta_{ij} m' M + \f_j^* \esti M \ge \esti \f_j^* M \ge 0 
\end{gather} 
where we write $A \ge B$ if for each irreducible $N$, the multiplicity 
 $[N:A]$ of $N$ in $A$ is not less than the multiplicity 
 $[N:B]$ of $N$ in $B$. 
Sum \eqref{gg} over $i,j \in \mu_q$, we get 
\begin{gather} 
a M + \Ind_{\Hlam{n-1}}^{\Hlam{n}} \Res_{\Hlam{n-1}}^{\Hlam{n}} M 
\ge \Res_{\Hlam{n}}^{\Hlam{n+1}} \Ind_{\Hlam{n}}^{\Hlam{n+1}} M 
\label{formula-a} 
\end{gather} 
for some $a \ge 0$. 
 
Next we claim that in $K(\Rlam)$ 
\begin{gather} 
(\sum \lambda_i) M + \Ind_{\Hlam{n-1}}^{\Hlam{n}} 
\Res_{\Hlam{n-1}}^{\Hlam{n}} M
= \Res_{\Hlam{n}}^{\Hlam{n+1}} \Ind_{\Hlam{n}}^{\Hlam{n+1}} M. 
\label{mack-cyc} 
\end{gather} 
Granting this for the moment, let $N$ be an irreducible $\Hlam{n}$-module 
with $N \neq M$. 
Then comparing the multiplicity of $N$ in \eqref{mack-cyc} and in 
\eqref{gg}, we see that all inequalities in \eqref{gg} must 
be equalities, and so the multiplicity of $N$ in  
$\f_j^* \esti M$ equals the multiplicity of $N$ in 
$\esti \f_j^* M$; so $[\esti, \f_j^*] (M)$ is a multiple of $M$. 
Furthermore, if $i \neq j$ then the multiplicity of $M$ in 
$\f_j^* \esti M$ and $\esti \f_j^* M$ is zero, as the central 
characters of $M$ and $\f_j^* \esti M$ differ. 
 
So to prove the proposition it remains to show \eqref{mack-cyc}. 
This is immediate from the Mackey formula for $\Hlam{n}$. 
A weak form of this may be immediately deduced from the 
Mackey formula for $\Haffn$; we omit further details. 
\end{proof} 
 
For $N$ an irreducible $\Hlam{n}$-module, write $\wti(N) = 
\nui(N) - \epsi(N)$, and 
$$ 
\wt(N) = \sum \wti(N) \Lambda_i \cc \: \mu_q \to \Z. 
$$       
Define a function $\delta_i\cc R^n \to \N$ by 
$$ \delta_i(s_1, \ldots, s_n) = \#\{a \mid 1 \le a \le n, s_a = i \},$$ 
and set 
$\wti(s) = -2 \delta_i(s) + \delta_{qi}(s) + \delta_{q^{-1} i}(s).$ 
Recall $\alpha_i = 2 \Lambda_i - \Lambda_{qi} - \Lambda_{q^{-1} i}$.

The following theorem which is a summary of the results of section 
\ref{sec-serre}, tells us that both $\epsi(N)$ and $\nui(N)$ may be read 
off the spectrum of $N$, and that their difference $\wti(N)$ depends 
only on the central character of $N$. 
\begin{thm} 
\label{thm-wt}
Let $N$ be an irreducible $\Hlam{n}$-module with central character $s$. Then 
\begin{enumerate} 
\item[(i)] $\wt(\ft N) = \wt(N) - \alpha_i$, if $\ft N \neq 0$. 
\item[(ii)] $\wti(\1) = \lambda_i$, where $\1$ is 
the irreducible $\Hlam{0}$-module. 
\item[(iii)] $\wti( N) =  \nui(N) - \epsi(N) = \lambda_i + \wti(s)$. 
\end{enumerate} 
\end{thm} 
\begin{proof} 
As every irreducible module is obtained from $\1$ by 
a sequence of raising operators $\widetilde{\f}_k$,  
it is clear that (i) and (ii) are equivalent to (iii). 
(ii) is immediate from the definition of $\fsti$ and the description 
$\Hlam{1} = R[x]/\prod (x-q^i)^{\lambda_i}$. 
So we prove (i). 
 
We first observe that as $N = \et \ft N$ if $\ft N \neq 0$, we 
have $\wti(\ft N) = \wti(N) - 2$. 
Further, as $[\e_j^*, \fsti] = 0 = [\f_j^*, \fsti]$ if 
$j \notin \{ qi, i, q^{-1}i \}$, 
for such $j \:$  $\wti(\widetilde{\f}_j N) = \wti(N)$, in agreement 
with (i). 
Hence for $q \neq -1$, the content of the theorem is the assertion
$$ 
        \wti(\widetilde{\f}_{qi}N) = \wti(N) + 1,\qquad  \text{ if }\quad 
\widetilde{\f}_{qi}N \neq 0. 
$$ 
This is immediate from proposition \ref{prop-handf}. 
Likewise, if $q = -1$, proposition \ref{prop-q-1} is equivalent 
to the assertion of the theorem. 
\end{proof} 
We are now in a position to finish the proofs of theorem \ref{thm-nu} 
and theorem \ref{thm-[e,f]}. 
For $M$ an irreducible module, let $\widetilde{\varphi}_i(M)$ be the integer 
defined in theorem \ref{thm-nu}. Then by theorem \ref{thm-nu} (iii),
\begin{align*}
\sum_{i \in \mu_q} \widetilde{\varphi}_i(M) = & \sum_{i,j \in \mu_q} \dim \Hom (\fsti M, 
f_j^* M) \\ 
= & \sum_{i,j\in \mu_q} \dim \Hom (M,\esti  f_j^* M) \\
= &  \dim \Hom(M,  \Res_{\Hlam{n}}^{\Hlam{n+1}} \Ind_{\Hlam{n}}^{\Hlam{n+1}} M)
\end{align*}
and by \eqref{mack-cyc}, this is 
$$ \leq \sum \lambda_i  + \sum \epsi(M). $$
But $\widetilde{\varphi}_i(M) \ge \nui(M)$ for all $i\in \mu_q$, by 
the last paragraph of  theorem \ref{thm-nu}, and 
by theorem \ref{thm-wt}{ (iii)}
$$ \sum  \lambda_i = \sum \left( \nui(M) - \epsi(M)\right)  .$$
It follows that
$$ \sum \lambda_i = \sum  \left( \nui(M) - \epsi(M)\right)
\leq \sum  \left( \widetilde{\varphi}(M) - \epsi(M)\right)
 \leq \sum  \lambda_i,
$$
and hence all the inequalities above are equalities. In particular
$\widetilde{\varphi}_i(M) = \nui(M)$, for all $i$, completing the proof
of theorem \ref{thm-nu}.

Finally, theorem \ref{thm-epsilon} (ii) and theorem \ref{thm-nu} (ii) now give
\begin{align*}
[ \esti\fsti M - \fsti\esti M : M]
= & (\epsi(M)+1)\nui(M) - (\nui(M) + 1)\epsi(M) \\
= & \nui(M) - \epsi(M)
\end{align*}
completing the proof of theorem \ref{thm-[e,f]}. 
%
 
To finish verifying the defining relations of $\sl$, we need 
only observe that theorem \ref{thm-[e,f]} immediately implies 
\begin{lemma} 
$[\hsti, \e_j^*] = c_{ij} \e_j^*, \qquad [\hsti, \f_j^*] = -c_{ij} \f_j^* \cc \: 
K(\Rlam) \to K(\Rlam)$, 
where $c_{ij} = 2 \delta_{ij} - \delta_{i, qj} - \delta_{qi, j}$ 
is the Cartan matrix of $\sl$. 
\end{lemma} 

\begin{rem} M.{} Vazirani has recently improved on theorem \ref{thm-[e,f]}
by determining the relation between $\esti\fsti M$ and $\fsti\esti M$
for $M$ irreducible, before passage to the Grothendieck group.
Her results extend and clarify theorems \ref{thm-nu}(iii) and \ref{thm-epsilon}(iii).
\end{rem}

\section{Uniqueness of the crystal}
\label{sec-unique}
In this section we determine the crystal graph of $K(\Raff)$.  This
admits many different combinatorial descriptions, each of which it is
possible to interpret Hecke-theoretically. Rather than do this, we
prove one more property of the crystal $B_\aff$ of $K(\Raff)$. This
property (proposition \ref{eps-hat}), together with what we have proved earlier, is
already sufficiently strong to show combinatorially the uniqueness of
the crystal. In fact other than \ref{eps-hat}, all we need is that the
crystals $\B$ admit a description purely in terms of the crystal
$B_\aff$ and the involution on $B_\aff$ induced by the
antiautomorphism $\sigma^*$.

Recall that we write $\ethat = \sigma^* \et \sigma^*$, and
also define $\esthat = \sigma^* \esti \sigma^*$
and $\fthat = \sigma^* \ft \sigma^*$.

\begin{prop}
\label{eps-hat}
Let $M$ be an irreducible $\Haff{M}$-module, and write $c = \epsihat(M)$
\begin{enumerate}
\item[(i)] Suppose $\epsihat(\ft M) = \epsihat(M)$. Then
$$ (\ethat)^c (\ft M) = \ft ( {\ethat}{}^c M) $$
\item[(ii)] If $\epsihat(\ft M) = \epsihat(M) + 1$, then $\ethat \ft M = M$.
\end{enumerate}
\end{prop}
\begin{proof}
(i) We have $M = (\fthat)^c N  = \cosoc \Ind (q^iK_c \boxtimes N)$,
where $N$ is an irreducible $\Haff{n}$ module with $\epsihat(N) = 0$.

Set $Q_a = (\esthat)^{c-a} \ft M$, so that in the Grothendieck
group $Q_a$ is some number of copies of $(\ethat)^{c-a} \ft M$ plus terms
with strictly smaller $\epsihat$. In particular we have that 
$\epsihat(A) \leq a$  for all $A$ that occur in $Q_a$, and
$Q_0$ is just some copies of $(\ethat)^c \ft M$.

We will show by decreasing induction on $a$ that there is a non-zero map
$$ \gamma_a : \Ind(q^iK_a \boxtimes N \boxtimes  q^iJ_1) \to Q_a.$$
If $a=c$, $Q_a = \ft M= \cosoc \Ind(M\boxtimes q^i J_1)$ is a
quotient of $\Ind(q^iK_c\boxtimes q^iJ_1)$ so our induction starts.
Now suppose $\gamma_a$ exists, and $a\geq 1$.
Consider $\Res^{a+n-1}_{1,a+n} \Ind (q^iK_a \boxtimes N\boxtimes q^i J_1)$.
By the Mackey formula, this has a three step filtration
$0 \subset F_1 \subset F_2 \subset F_3$ with successive quotients
\begin{multline*}
 F_1 = \Ind^{1,a+n}_{1,a-1,n,1} \Res^{a,n,1}_{1,a-1,n,1} (q^iK_a \boxtimes N 
\boxtimes q^iJ_1), \\
F_2/F_1 = \Ind^{1,a+n}_{1,a,n-1,1} {}^w \Res^{a,n,1}_{a,1,n-1,1} (q^iK_a \boxtimes N 
\boxtimes q^iJ_1), \\
F_3/F_2 = \Ind^{1,a+n}_{1,a,n}(q^iJ_1 \boxtimes q^i K_a \boxtimes N),
\end{multline*}
where $w$ is the obvious permutation.

As $\gamma_a \neq 0$, Frobenius reciprocity gives a copy of 
$q^i K_a \boxtimes N \boxtimes q^iJ_1 $ in the image of $\gamma_a$, and so
$$\tilde{\gamma}_a := \esthat \gamma_a : \esthat \Ind (q^i K_a \boxtimes N
\boxtimes q^iJ_1) \to \esthat Q_a = Q_{a-1} $$ 
is non-zero. Suppose that $\tilde{\gamma}_a$ is zero when restricted to the 
$q^i$-eigenspace of $X_1$ on $F_1$. 
As there is no $q^i$-eigenspace of $X_1$ on $F_2/F_1$, we must have
a non-zero homomorphism from $F_3/F_2$ to $Q_{a-1}$, i.e.{} a non-zero
homomorphism
$$ \Ind(q^iK_a\boxtimes N) \to Q_{a-1}.$$
But $\epsihat(\cosoc \Ind(q^iK_a\boxtimes N)) = a > \epsihat(A)$,
for any constituent $A$ of $Q_{a-1}$. So this is not possible, and it
must be that $\tilde{\gamma}_a$ restricts to a non-zero homomorphism
on the $q^i$-eigenspace of $X_1$ on $F_1$.
As $\esthat (q^iK_a)$ has a filtration 
with subquotients $q^iK_{a-1}$, there must be a non-zero map 
$\gamma_{a-1} : \Ind(q^iK_{a-1} \boxtimes N \boxtimes q^iJ_1)
\to Q_{a-1}$.

We now take $a=0$ and conclude there is a non-zero homomorphism
$$ \gamma_0 : \Ind(N \boxtimes q^i J_1) \to Q_0, $$
hence $\ft (\ethat)^c M = \ft N = \cosoc \Ind (N \boxtimes q^iJ_1)$
is a subquotient of $Q_0$. But $Q_0$ is a multiple of $(\ethat)^c\ft M$,
so we have indeed shown that $\ft (\ethat)^c M = (\ethat)^c\ft M$.

(ii) Again write $N = (\fthat)^c M$.  As multiplication is commutative
in the bialgebra $K(\Raff)$, $\Ind(q^iK_c \boxtimes N \boxtimes
q^iJ_1)$ equals $\Ind(q^iK_c \boxtimes q^iJ_1 \boxtimes N) =
\Ind(q^iK_{c+1} \boxtimes N)$ in the Grothendieck group. 
Hence $\Ind(q^iK_c \boxtimes N \boxtimes q^iJ_1)$
is $(\fthat)^{c+1} N = \fthat M$ plus terms $A$ with $\epsihat(A) \leq c$.
As  $\Ind(q^iK_c \boxtimes N \boxtimes
q^iJ_1)$ surjects onto $\ft M$, if $\epsihat (\ft M) = c+1$ it
must be that $\fthat M = \ft M$. 
\end{proof}

This is the last property we will need of the representations of 
the affine Hecke algebra. The next proposition is a formal consequence of
what we have already proved.

For $M$ an irreducible $\Haff{n}$-module with central character $s$,
define $$ \wt_i'(M) = \wt_i(s).$$
Also define $(\ethat)^{\max}(M) = (\ethat)^c(M)$, where $c = \epsihat(M)$.

\begin{prop}
\label{eps-crys}
Let $M$ be an irreducible $\Haff{m}$-module, and
 write $c = \epsihat(M)$, $\bar{M} = (\ethat)^{\max}(M)$.
\begin{enumerate}
\item[(i)] $\epsi(M) = \max(\epsi(\bar{M}), c - \wt_i'(\bar{M}))$.
\item[(ii)] Suppose $\epsi(M) > 0$. Then 
$$\epsihat(\et M) = \begin{cases} c, & \text{ if } \epsi(\bar{M}) \geq c - \wt_i'(\bar{M}) \\
c-1, & \text{ if } \epsi(\bar{M}) < c - \wt_i'(\bar{M}) \end{cases} 
$$
\item[(iii)] Suppose $\epsi(M) > 0$. Then
$$(\ethat)^{\max}(\et M) = \begin{cases} \et(\bar{M}),  & \text{ if } \epsi(\bar{M}) \geq c - \wt_i'(\bar{M}) \\
\bar{M}, &  \text{ if } \epsi(\bar{M}) < c - \wt_i'(\bar{M}) \end{cases}
$$
\end{enumerate}
\end{prop}
\begin{proof}
Let $\lambda : \mu_q \to \N$ be such that $\sum\lambda_i < \infty$,
and suppose $N \in \Rep \Hlam{a}$. Then corollary \ref{cor-epscheck} tells us that
$k = \nui(N)$ means 
$$ \epsihat(N) \leq \epsihat(\ft^k N) = \lambda_i, \qquad \text{and }
\epsihat(\ft^{k+1}N) = \lambda_i + 1,$$
and that theorem \ref{thm-wt} shows that $\nui(N) = \lambda_i + \epsi(N) +
\wt_i'(N)$.

Take $N = \et^m M$, where $m = \epsi(M)$, and define 
$\lambda(y) : \mu_q \to \N$ by setting $\lambda_i(y) = \epsihat(N) + y$,
and setting $\lambda_j(a)$ to be any integer much greater than $m + y + a$,
when $j \neq i$. It then follows from the previous paragraph applied
to $\lambda(0), \lambda(1), \dots$ that for all $s$
$$ \epsihat(\ft^s N) = \begin{cases} 
\epsihat(N), \qquad &  s \leq \epsihat(N) + \epsi(N) + \wt_i'(N) \\
s - \epsi(N) - \wt_i'(N), \qquad & s \geq \epsihat(N) + \epsi(N) + \wt_i'(N) 
\end{cases}
$$
and as $\epsi(N) = 0$ it follows that 
$$ \epsihat(\ft^m N) = \max(\epsihat(N), m - \wt_i'(N)), $$
i.e.{} that $\epsihat(M) = \max(\epsihat(\et^{\epsi(M)}M), 
\epsi(M) - \wt_i'(\et^{\epsi(M)}M))$. Applying this to $\sigma^* M$ we get (i).

To see (ii), observe that $\epsihat(\et M) = \epsihat(M) - 1 $
precisely when 
$$ m = \epsi(M) > \wt_i'(N) + \epsihat (N) $$
and that otherwise $\epsihat(\et M)= \epsihat(M)$. But $\wt_i'(N) = \wt_i'(M)+ 2m$, and so (ii) follows if we show that 
$\wt_i'(M) + \epsihat(N) + m < 0$ if and only if $wt_i'(M) +  \epsi(\bar{M}) + c < 0$. But $\wt_i'(M) + \epsihat(N) + m = \max(\wt_i'(M) + \epsihat(N) + \epsi(\bar{M}),\epsihat(N) -c)$, by (i), and $\epsihat(N) - c \leq 0$ always.
Similarly $wt_i'(M) +  \epsi(\bar{M}) + c = \max(\wt_i'(M) + \epsihat(N) + \epsi(\bar{M}),\epsi(\bar{M}) - m)$, and $\epsi(\bar{M}) - m\leq 0$ always. 
Finally, proposition \ref{eps-hat}i shows that $\epsi( \bar{M}) =m $ if and only if
$ \epsihat(N) =c$; hence (ii).

Now (iii) is immediate from (ii) and proposition \ref{eps-hat}.
\end{proof}

\subsection{Combinatorial consequences} 
We now show that we have enough properties to completely describe the
tensor category of ``integrable lowest weight crystals'', and hence to
describe the crystals themselves. We follow \cite{Ks}, especially 8.2.
(This is a kind of purely combinatorial Tannakian
property.)

So we recall some ideas of Kashiwara \cite{Ks}. Recall that if $B_1$
and $B_2$ are two crystals, their tensor product $B_1 \otimes B_2$ is
the crystal whose underlying set is $B_1 \times B_2$ and with
$$ \et(b_1 \otimes b_2) = \begin{cases} \et b_1 \otimes b_2 &
\quad \text{ if } \nui(b_1) \geq \epsi(b_2) \\
b_1 \otimes \et b_2 &
\quad \text{ if } \nui(b_1) < \epsi(b_2) \end{cases} $$
$$ \ft(b_1 \otimes b_2) = \begin{cases} \ft b_1 \otimes b_2 &
\quad \text{ if } \nui(b_1) >\epsi(b_2) \\
b_1 \otimes \ft b_2 &
\quad \text{ if } \nui(b_1) \leq \epsi(b_2) \end{cases} 
$$

Define the crystal $B_i$, for $i \in \mu_q$, to have underlying set
$\{ b_i(n) \mid n \in \Z \}$ and set
$$ \varepsilon_j(b_i(n)) = \begin{cases} - n, & j = i \\ -\infty & j \neq i,  \end{cases} \qquad
 \varphi_j(b_i(n)) = \begin{cases}  n, & j = i \\ -\infty & j \neq i,  \end{cases} 
$$
$$
\widetilde{\e}_j( b_i(n)) = \begin{cases} b_i(n+1), & j = i \\ 0 & j \neq i, \end{cases}
\qquad
\widetilde{\f}_j( b_i(n)) = \begin{cases} b_i(n-1), & j = i \\ 0 & j \neq i \end{cases}
$$
and write $b_i = b_i(0)$. 

Recall that a strict embedding of crystals is an injective map
$\psi:B_1 \inj B_2$ such that $\psi$ commutes with $\et$ and $\ft$
for all $i \in \mu_q$. Define $B_\infty$ to be the same crystal as 
$B_\aff$, except that we set $\nui(M) := \epsi(M) + \wt_i(s)$ if
$M$ has central character $s$.

We can now rephrase proposition \ref{eps-crys} as

\begin{prop} For each $i \in \mu_q$, define a map $\Psi_i :  B_\infty \to 
B_\infty \otimes B_i$ by sending $M$ to $(\esthat)^c(M) \otimes \ft^c b_i$,
where $c = \epsihat(M)$. Then $\Psi_i$ is a strict embedding of crystals.
\end{prop}

It is a result of Kashiwara that this determines the crystal
$B_\infty$ (see \cite{Ks} and proposition 3.2.3 of \cite{KS}). For the
convenience of the reader we reproduce the argument. Choose a sequence
$(i_1,i_2, \dots)$ in $\mu_q$ so that each $i\in\mu_q$ appears
infinitely often.  Define a map $\Phi_n : B_\infty \to B_\infty
\otimes B_{i_n} \otimes \cdots \otimes B_{i_1}$ by $\Phi_n =
\Psi_{i_n}\circ \cdots \circ \Psi_{i_1}$.  Then for any $b \in
B_\infty$ there exists an $n$ such that $\Phi_n(b) = b_0 \otimes
\widetilde{\f}^{a_n}_{i_n}b_{i_n}\otimes\cdots\otimes
\widetilde{\f}^{a_1}_{i_1}b_{i_1}$. The sequence $(a_1,a_2,
\cdots,a_n,0,0,\cdots)$ does not depend on $n$. This embeds $B_\infty$
as the smallest subcrystal of $B_{\text{Kas}}$ containing
$(0,0,\cdots)$, where we define $B_{\text{Kas}}$ to be the crystal
whose underlying set is the set of sequences $\{(a_i) \in \Z \mid a_i
= 0 \text{ for } i \gg 0 \}$, and whose crystal structure is defined
by sending $(a_1,\dots,a_n,0,\dots)$ to $\cdots \otimes
b_{i_2}(-a_{i_2})\otimes b_{i_1}(-a_{i_1})$.  Hence $ B_\infty$ is
completely determined by the Dynkin diagram $\mu_q$, as
$B_{\text{Kas}}$ is.

\section{Reaping the harvest}
\label{sec-reap}
In this section we summarise the theorems of the previous sections and
identify the Hopf algebra $K(\Raff)$ and its comodules $K(\Rlam)$.  As
a consequence of this rigid structure, we obtain a parameterization of
irreducible modules for $\Haff{n}$ and $H_n^\lambda$ over any field
$R$ that depends only on $n$, $\lambda$ and $\ell = |\mu_q|$ (and not
on $R$, or even the characteristic of $R$).  We observe that the
crystal graph of $B_\lambda$ we have defined coincides with the
crystal graph defined by Kashiwara and Lusztig.   We emphasise that there
is no further Hecke theoretic content in the parameterisation of
representations---any of the many combinatorial descriptions of the
crystal basis gives a combinatorial parameterisation of
$B_\lambda$. Thus we can label modules by tuples of partitions, or
Littelmann paths, or paths in a perfect crystal, for example. The
identification of this with Deligne-Langlands parameters \cite{KL,G}
is a pleasant exercise.


For convenience we dualize $K(\Rlam)$, and denote the
adjoints to $\edst$, $\fdst$, $\binom{h_i^*}{n}$ by 
$\ed$, $\fd$, $\binom{h_i}{n}$. Also write $\1$ for the generator of 
$K(\Rep H_0^\lambda) = \Z$ dual to the trivial representation
of $H_0^\lambda$.
\begin{thm}
\label{uz-bi}
The map $\U_\Z \n_\ell \to K(\Raff)^\ast$ which sends the
generators $\ed$ to the elements $\ed {\bf 1}_{\aff} = \delta_{q^i K_n}$ is an isomorphism
of bialgebras.
\end{thm}
\begin{thm}
\label{uz-hw}
The operators $\ed$, $\fd$, $\binom{h_i}{n}$ define a structure of a
$\U_\Z \sl$-module on $K(\Rlam)^\ast$.
The module $K(\Rlam)^\ast$ is a $\Z$-form of the
irreducible integrable lowest weight module for $\sl$ with
lowest weight $\lambda$ (and lowest weight vector $\1$).
Under this identification, the Shapovalov form on the 
module becomes the form of section \ref{sec-shap}.
\end{thm}
\begin{rem} The bialgebra structure on $K(\Rfin)^*$
is that given by the principal realisation 
of the basic representation, i.e.~ the identification of this with
the Hopf algebra $\Z[x_i \mid \ell \nshortmid i]$.
\end{rem}

\begin{thm}
\label{thm-crystal}
The crystal graph of the $\sl$-module $K(\Rlam)^\ast$
is the graph $(\B, \e_i, \f_i)$ defined in section \ref{sec-crystal}.
\end{thm}
We prove theorems \ref{uz-bi} and \ref{uz-hw} simultaneously.

Step 1:
As $K(\Rlam)^\ast$ is a torsion-free $\Z$-module, and 
$\ed, \fd, \binom{h_i}{n}$ are well defined operators which
are determined by the actions of $e_i$, $f_i$ and $h_i$,
it is enough to
check that the defining relations on $\e_i$, $\f_i$, $\h_i$
are satisfied.
This was done in the last section, and so we have a well defined
action of $\U_\Z \sl$ on $K(\Rlam)^\ast$.
Similarly, we have  a well defined bialgebra morphism
$\U_\Z \n_\ell \to K(\Raff)^\ast$.

Step 2: 
The map $\U_\Z \n_\ell \to K(\Raff)^\ast$ is surjective.
This follows from lemma \ref{lem-shap3}, which is precisely the statement that 
the map $\U_\Z \n_\ell \to K(\Rlam)^*$ is surjective, and the fact that
every irreducible $\Haffn$-module is a $\Hlam{n}$ module for 
some $\lambda$.

Step 3: We show $K(\Rlam)_\Q^\ast$ is the irreducible $\sl$-module
with lowest weight $\lambda$. But $K(\Rlam)_\Q^\ast$ is a
$\U_\Q\sl$-module on which $e_i$ and $f_i$ act locally nilpotently
(theorems \ref{thm-epsilon} and \ref{thm-nu}), and on which $h_i$ acts
semisimply. By step 2 it is generated by $\1$ as a
$\U_\Z\n_\ell$-module. Hence it is an irreducible integrable lowest
weight module. By theorem \ref{thm-wt} (ii) this weight is $\lambda$.

Step 4:
Finally we show the map $\U_\Z \n_\ell \to K(\Raff)^\ast$ is injective.
Let $x$ be in its kernel. As the action of $\U_\Z\n_\ell$ on
$K(\Rlam)^*$ factors through $K(\Raff)^*$, $x$ acts as zero on every 
$K(\Rlam)^*$ hence on every integrable lowest weight module. It follows 
that $x = 0$.

This concludes the proof of theorems \ref{uz-bi} and \ref{uz-hw}.
Theorem \ref{thm-crystal} follows from the description of the crystal
$B_\aff$ in section \ref{sec-unique}, and the fact that the crystal
of $\U\n$ has the same description \cite{Ks}.

\subsection{The $p$-canonical basis}
The above theorems prompt the following definition. 
Fix a prime $p \ge 0$, and a positive integer $\ell$. 
\begin{defn} {\it The $p$-canonical basis\/} of $\sl$ 
is the basis of the  module $\U_\Z { \n}_\ell$ given 
by the dual of the irreducible $\Haffn$-modules, where 
$\Haffn$ is the affine Hecke algebra over $R = \overline{\F}_p$, 
an algebraically closed field of characteristic $p$, 
and $q \in R$ is a primitive $\ell^{th}$ root of unity. 
(If $p=0$, take $R=\C$.) 
\end{defn} 
 
(If $(p,\ell) \neq 1$, we can still define such a basis. 
The case $p=\ell$ is explained in the next section.) 
This basis has the following pleasant properties, among 
others. 
\begin{enumerate} 
\item[(i)] The basis of $\U_\Z{\n}_\ell$ descends 
to give a basis for all integrable lowest weight modules. 
\item[(ii)] The structure constants of $\e_i$ and $\f_i$ on 
this basis are non-negative integers. 
\item[(iii)] The $0$-canonical basis is a non-negative integral 
combination of the $p$-canonical basis, for each prime $p$. 
\item[(iv)] The  $0$-canonical basis is  the canonical basis 
($=$ {\it global crystal basis\/}) of Lusztig and Kashiwara. 
\end{enumerate}

Note that we have given (elementary!) proofs of (i)--(iii) 
in this paper; property (iv) is immediate from the 
Kazhdan-Lusztig description of the affine Hecke algebra in 
geometric terms 
and Lusztig's definition of the canonical basis in terms of 
perverse sheaves on quivers; this is explained 
(tersely) in \cite{G}.
\section{Modifications when $q=1$.} 
\label{sec-q=1}
If $q=1$ all of the above theorems and constructions go through, 
without change, once we make the appropriate definitions. 
 
Define $\mu_q$ to be the image of $\Z \to R$, $1 \mapsto 1$, so 
$\ell$ is the characteristic of $R$, $\ell \in \N \cup \{\infty\}$. 
Instead of the affine Hecke algebra, we work with the 
{\it degenerate\/} or {\it graded\/} affine Hecke algebra  
$\Hdeg$, defined by Drinfeld \cite{D} and Lusztig \cite{L}. 
This is isomorphic as an $R$-module to 
$$ 
        R[S_n] \otimes R[X_1, \ldots X_n] 
$$ 
with algebra structure defined by requiring that 
$R[S_n]$ and $R[X_i]$ are subalgebras, and that 
$$ 
        s_i\cdot f -  {}^{s_i}\!f \cdot s_i =  
        \frac{f - {}^{s_i}\!f}  { X_i -X_{i+1}}. 
$$ 
Given a function  $\lambda : \mu_q \to \Z_+$  
such that  
$\sum \lambda_i < \infty$, define the  
{\it degenerate cyclotomic algebra\/} $\Hldeg = \Hdeg/I_\lambda$, 
where $I_\lambda = \Hdeg \cdot \prod_{i \in \mu_q}(X_1 - i)^{\lambda_i} \cdot
\Hdeg$. 
If $\lambda = \Lambda_0$, $\overline{H}_n^{\Lambda_0} = R[S_n]$; 
in general one shows that $\dim_R \Hldeg = r^n n!$, 
where $r = \sum \lambda_i$. 
We define $\Rep_q \Hdeg$ to be the subcategory of  
$\Hdeg$-modules on which $X_1$ acts with eigenvalues in $\mu_q$. 
All other definitions and theorems are as before, once we agree 
to write the group law in $\Im \Z = \mu_q$ multiplicatively, 
so $qi$ denotes what is usually written $i+1$. 
No changes are necessary in the proofs. 
 
In particular, take $R = \F_p$, to see that 
$\widehat{\mathfrak sl}_p$ controls the representation 
theory of the symmetric group in characteristic $p$.


 

\end{document}